\documentclass[a4paper,11pt]{amsart}
\usepackage{amssymb,amsfonts,amsxtra}
\usepackage[all]{xy}
\usepackage{fullpage}
\newtheorem{theorem}{Theorem}[section]
\newtheorem{cor}[theorem]{Corollary}

\newtheorem{prop}[theorem]{Proposition}

\theoremstyle{definition}
\newtheorem{defi}[theorem]{Definition}
\newtheorem{rem}[theorem]{Remark}
\numberwithin{equation}{section} 
\thanks{The author was partially supported by the EPSRC grant No. GR/R84276/01}
\begin{document}
\title[Hochschild cohomology and moduli spaces...]
{Hochschild cohomology and moduli spaces of strongly homotopy
associative algebras}
\author{A.Lazarev }
\address{Mathematics Department, University of Bristol, Bristol, BS8 1TW,
England.} \email{a.lazarev@bristol.ac.uk}
\begin{abstract}Motivated by ideas from stable homotopy theory we
study the space of strongly homotopy associative multiplications
on a two-cell chain complex. In the simplest case this moduli
space is isomorphic to the set of orbits of a group of invertible
power series acting on a certain space. The Hochschild cohomology
rings of resulting $A_\infty$-algebras have an interpretation as
totally ramified extensions of discrete valuation rings. All
$A_\infty$-algebras are supposed to be unital and we give a
detailed analysis of unital structures which is of independent
interest.\end{abstract} \maketitle
\section{Introduction} The notion of a strongly homotopy associative algebra
or of an $A_\infty$-algebra was introduced in  \cite{Sta} by
Stasheff and was recently much studied in connection with
deformation quantization and the Deligne conjecture, cf.
\cite{Kon1}, \cite{MS}. From the point of view of homotopy theory,
an $A_\infty$-algebra is the same as a differential graded algebra
(dga). However, for the purposes of explicit computations, it is
often more convenient to work with $A_\infty$-algebras rather than
with dga's.

The purpose of this paper is to study  `homotopy invariant' or
`derived' moduli spaces for $A_\infty$-algebras. It should be
noted that other authors also considered the problem of
constructing derived moduli spaces. Here we mention the works of
M.Schlessinger and J.Stasheff, cf. \cite{SS} and of V.Hinich,
\cite{Hin1},\cite{Hin}. Another approach making heavy use of
simplicial methods and homotopical algebra is developed in
\cite{BDG}. The case of $A_\infty$-algebras considered here,
exhibits, on the one hand, most of the representative features of
derived moduli space theory and, on the other hand, allows one to
perform concrete computations without the need of too much
apparatus.

Though our examples are purely algebraic, they are motivated by
the study of complex-oriented cohomology theories. There is a
parallel notion of an $S$-algebra, or an $A_\infty$-ring spectrum
in stable homotopy theory, cf. \cite{EKMM}. Many important spectra
of algebraic topology, especially those related to complex
cobordisms admit structures of $A_\infty$-ring spectra. The
existence of such structures was proved in  \cite{Laz1},
\cite{Rob} and in \cite{Goe} by methods of obstruction theory, but
up until now there was no attempt to classify such structures. Of
particular interest is the space of $A_\infty$ structures on
$K(n)$'s, the higher Morava $K$-theories. A closely related
problem is computation of $THH^\ast(K(n),K(n))$, the topological
Hochschild cohomology of $K(n)$. It should be noted that in the
topological situation one is forced to work in the abstract
setting of Quillen closed model categories and this makes the
classification problem much more difficult then the corresponding
algebraic one. Therefore it seems natural to consider the
algebraic problem first and this is what we do in the present
paper.

We now outline the problem under consideration and our approach.
Consider the field of $p$ elements $\mathbb{F}_p$ as an algebra
over $\mathbb{Z}$. Then we could form the derived Hochschild
cohomology of $\mathbb{F}_p$ with coefficients in itself as a
$\mathbb{Z}$-algebra (sometimes called Shukla cohomology, cf.
\cite{Shu}): $$HH_{\mathbb{Z}}(\mathbb{F}_p,\mathbb{F}_p):=
RHom_{\mathbb{F}_p\otimes^L\mathbb{F}_p}(\mathbb{F}_p,\mathbb{F}_p).$$
Here $\otimes^L$ denotes the derived tensor product over
$\mathbb{Z}$ and $RHom$ denotes the derived module of
homomorphisms. An easy computation then shows that
$HH^\ast_{\mathbb{Z}}(\mathbb{F}_p,\mathbb{F}_p)=\mathbb{F}_p[[z]]$
with $z$ having cohomological degree $2$.

This result is valid because $\mathbb{F}_p$ has a unique structure
of a $\mathbb{Z}$-algebra which happens to be commutative. If $R$
is an evenly graded commutative ring and $x\in R$ is a homogeneous
element that is not a zero divisor we obtain similarly
$HH_R^\ast(R/x,R/x)=R/x[[z]]$ with $z$ having cohomological degree
$|x|+2$. Notice, however, that we implicitly resolved the
$R$-algebra $R/x$ by the differential graded algebra
$\Lambda_R(y)$, with one generator $y$ in degree $|x|+1$ whose
square is $0$ and $dy=x$. In other words we  assumed that $R/x$ is
given the usual structure of an $R$-algebra, in particular that it
is commutative. In general there are many different structures of
an $A_\infty$-algebra on the complex $R\stackrel{x}{\rightarrow}
R$ (which is a model for $R/x$ in the derived category of
$R$-modules). Take for instance $R=\mathbb{Z}[v,v^{-1}]$ where the
element $v$ has degree $2$ and $x$ is a prime number $p\neq 2$.
Then the differential graded algebra $R[y]/(y^2-v)$ with
differential $dy=p$ is a model for $R/x$ (that is, its homology
ring is $R/x$), but it is not commutative even up to homotopy. It
turns out that its Hochschild cohomology is
$\hat{\mathbb{Z}}_p[v,v^{-1}]$ where $\hat{\mathbb{Z}}_p$ is the
ring of $p$-adic integers. In general for any structure of an
$A_\infty$-$R$-algebra on $R/p$ the ring $HH_R^\ast(R/p,R/p)$ is
filtered and complete with associated graded isomorphic to the
formal power series algebra over $\mathbb{F}_p[v,v^{-1}]$. It
follows that the ring $HH_R^\ast(R/p,R/p)$ is either an
$\mathbb{F}_p$-algebra or it has no $p$-torsion. The torsion-free
case corresponds to totally ramified extensions of the field of
the $p$-adic numbers (cf.\cite{Ser}) and we show that by varying
$A_\infty$ structures one can obtain extensions of arbitrary
ramification index that is coprime to $p$, the so-called tamely
ramified extensions.

So the natural problem is now to classify all possible $A_\infty$
structures of $R/x$ (or, equivalently all differential graded
$R$-algebras whose homology ring is $R/x$). We consider this as
part of a more general problem, namely the classification of all
$A_\infty$ structures on a two cell complex
$\{\Sigma^dR\stackrel{\partial}{\rightarrow} R \}$ with no
restrictions on $d$ or the differential $\partial$. An
$A_\infty$-algebra of this sort is called a Moore algebra as being
analogous to the Moore spectrum of stable homotopy theory.  The
resulting theory bears striking resemblance with the theory of
one-dimensional formal groups, cf. \cite{Haz} although we could
not establish a direct link between the two. There are essentially
two different cases: $d$ odd and $d$ even. We concentrate mainly
on the even case since it is most relevant to the parallel
topological problem. However, the odd case is also quite
interesting and apparently related to the homology of moduli
spaces of algebraic curves, cf. \cite{Kon2}.

We note that there is a close connection between the results of
this paper and algebraic deformation theory, cf. \cite{Ge},
\cite{PS}. For example Theorem \ref{norm} implies that any
deformation of a unital $A_\infty$-algebra is equivalent to a
unital one. This and other implications of our results in
deformation theory will be explained elsewhere.

We need to make a comment about grading. Throughout the paper we
work with $\mathbb{Z}$-graded complexes of modules over a
$\mathbb{Z}$-graded even commutative ring $R$. Some other authors,
e.g. \cite{Kon} work in the slightly less general
$\mathbb{Z}/2\mathbb{Z}$-graded context. These two approaches are
closely related. If we have a $\mathbb{Z}$-graded object we could
always forget down to a $\mathbb{Z}/2\mathbb{Z}$-graded object.
Conversely, tensoring everything in sight with the ring
$\mathbb{Z}[v,v^{-1}]$ with $|v|=2$ we obtain a
$\mathbb{Z}$-graded object from a $\mathbb{Z}/2\mathbb{Z}$-graded
one.  This procedure is routinely employed in topology when
studying complex oriented cohomology theories, cf. for example
\cite{AHS}.

This paper is organized as follows. In section $2$ we recall the
definition of an $A_\infty$-algebra and $A_\infty$-morphism and
collect various formulae which will be needed later on. The
material presented here is fairly standard, except that we define
$A_\infty$-algebras over commutative graded rings rather than over
fields as normally done.

In section $3$ we study unital structures and prove a formula for
the action of a unital automorphism on a given $A_\infty$
structure. It turns out that this and other formulae in the theory
of $A_\infty$-algebras are best handled using the language of the
dual cobar-construction (or, perhaps,  of the `cobar-construction
of the dual') which could be thought of as a formal noncommutative
(super)manifold in the sense of \cite{Kon2} . The seemingly
trivial passage from bar to cobar construction is the main
technical invention of this paper. We hope that it will have
further applications.

In section $4$ we define the Hochschild complex for
$A_\infty$-algebras and prove that it is  homotopically equivalent
to a normalized complex. This theorem is well-known for (strictly)
associative algebras.

Section $5$ introduces  Moore algebras which are our main object
of study. Moore algebras are in several respects similar to
one-dimensional formal group laws and we prove the analogue of
Lazard's theorem stating that the functor associating to a ring
the set of Moore algebras over it is representable by a certain
polynomial algebra on infinitely many generators.

In section $6$ we consider the problem of classification of even
Moore algebras over a field or a complete discrete valuation ring.
This problem is equivalent to the classification of orbits of a
certain action of the group of formal power series without the
constant term. We obtain complete classification in characteristic
zero and some partial results in characteristic $p$.

In section $7$ we compute Hochschild cohomology of even Moore
algebras and discuss its relation with totally ramified extensions
of discrete valuation rings.

{\bf Notations and conventions.} In sections $2-4$ we work over a
fixed evenly graded commutative ground ring $R$. The symbols $Hom$
and $\otimes$ always mean $Hom_R$ and $\otimes_R$. In Sections
$5-7$ the emphasis is shifted somewhat in that the ground ring is
varied. It is still denoted by $R$ sometimes with subscripts, e.g.
$R_e$ and $R_o$, and we use unadorned $Hom$ and $\otimes$ where it
does not cause confusion.

A graded ring whose homogeneous nonzero elements are invertible is
referred to as a graded field. A graded discrete valuation ring is
a Noetherean graded ring having a unique homogeneous ideal
generated by a nonnilpotent element.

The set of invertible elements in a ring $R$ is denoted by
$R^\times$.

{\it Acknowledgement.} The author wishes to thank J.Greenlees
whose visit to Bristol in February $2001$ provided stimulus for
the author to start this project and J. Rickard and J. Chuang for
many useful discussions.

\section{Preliminaries}
In this section we collect the necessary definitions and facts
about $A_\infty$-algebras. The details may be found in the
definitive monograph by Markl, Shnider and Stasheff \cite{MSS} or
in the papers by Getzler-Jones \cite{GJ}, and Keller \cite{Kel}.

 Let $R$ be an evenly graded commutative ring.  Since we need to
work in the derived category of unbounded complexes of $R$-modules
we will recall some basic facts following \cite{KM}. An $n$-sphere
$R$-module is a free $R$-module on one generator in degree $n$. A
cell $R$-module is the union of an expanding sequence of
$R$-submodules $M_n$ such that $M_0=0$ and $M_{n+1}$ is the
mapping cone of a map $\phi_n:F_n\rightarrow M_n$ where $F_n$ is a
direct sum of sphere modules (perhaps of different degrees). Thus
a cell $R$-module is necessarily a complex of free $R$-modules.
Conversely, it is easy to see that a bounded below complex of free
$R$-modules is a cell $R$-module. In our applications we will be
concerned with only such cell complexes.

Further, any $R$-module $M$ admits a cell approximation, that is
there is a cell $R$-module $\Gamma{M}$ and a a quasi-isomorphism
of complexes $\Gamma{M}\rightarrow M$. The functors $?\otimes M$
and $Hom(M,?)$ where $M$ is a cell $R$-module preserve
quasi-isomorphisms and exact sequences in the variable $?$. This
allows one to define the derived functors of $\otimes$ and $Hom$
by setting $M\otimes^LN:=\Gamma M\otimes N$ and
$RHom(M,N):=Hom(\Gamma M,N)$.

Let $A$ be a cell $R$-module and ${T}A$ the tensor algebra of
$A$:$${T}A=R\oplus A\oplus A^{\otimes 2}\oplus\ldots.$$ Then
${T}A$ is a coalgebra via the comultiplication
$\Delta:{T}A\rightarrow {T}A\otimes {T}A$ where
$$\Delta(a_1,\ldots,a_n)=\sum (a_1,\ldots, a_i)\otimes
(a_{i+1},\ldots,a_n).$$ It is standard that a coderivation
$\xi:TA\rightarrow TA$ of the coalgebra $TA$ is determined by the
composition $TA\stackrel{\xi}{\rightarrow} TA\rightarrow A$ where
the second map is the canonical projection. Denoting the
components of this composite map by $\xi_i$ we see that a
coderivation $\xi\in Coder(TA)$ is specified by a collection of
maps $\xi_i:A^{\otimes i}\rightarrow A$.
 Let us introduce a
filtration on the coalgebra $TA$ by setting
$F_pTA=\oplus_{i=0}^{p}T^iA$. Then the space of coderivations of
$TA$ also acquires a filtration. Namely, we say that a
coderivation $\xi$ has weight $\geq p$ if $\xi(F_iTA)\subset
F_{i-p}TA$. Then clearly the elements of weight $\geq -1$ form the
whole space $Coder(TA)$ and filtration by weight is an exhaustive
Hausdorff filtration on the space $Coder(TA)$. We will denote the
set of coderivations of weight $\geq p$ by $O(p)$.

For a graded $R$-module $A$ we will denote by $\Sigma A$ its
suspension: $(\Sigma A)_i=A_{i-1}$.
\begin{defi}The structure of an $A_\infty$-algebra on a cell
 $R$-module $A$ is a coderivation $m:{T}\Sigma A\rightarrow {T}\Sigma A$ of
degree $-1$ such that $m^2=0$, $m(T^0\Sigma A)=0$ and the first
component $m_1$ of $m$ is (the suspension of) the original
differential on $A$. Thus an $A_\infty$-algebra is a pair $(A,m)$.
We will frequently omit mentioning $m$  and simply refer to the
$A_\infty$-algebra $A$.
\end{defi}
\begin{rem} The condition that $m(T^0\Sigma A)=0$ means that $m\in
O(0)$, or that the zeroth component $m_0$ of $m$ vanishes. Some
authors consider $A_\infty$-algebras with nonvanishing $m_0$, cf.
\cite{GJ1}.\end{rem}This definition is slightly more general then
the usual one in that $R$ is not assumed to be a field. We
emphasize here that the grading on the ground ring $R$ will be
essential for our constructions. For an $A_\infty$-algebra $A$ we
will call the coalgebra ${T}\Sigma A$ with the differential $m$
the bar-construction of $A$ and use the symbol $BA$ to denote it.
Following the usual tradition we will denote the element $\Sigma
a_1\otimes \Sigma a_2\ldots\Sigma a_n\in (\Sigma A)^{\otimes n}$
by $[a_1|a_2|\ldots |a_n]$. The following formula shows how to
recover the coderivation $m:T\Sigma A\rightarrow T\Sigma A$ from
its components:
\begin{equation}\label{rec1}m[a_1|\ldots|a_{n}]=\sum_{k=1}^{n}\sum_{i=0}^{n-k}
(-1)^{|a_1|+\ldots+|a_i|+i}[a_1|\ldots|a_i|
m_k[a_{i+1}|\ldots|a_{i+k}]|a_{i+k+1}|\ldots|a_{n}]\end{equation}

The components ${m_i}:\Sigma A^{\otimes i}\rightarrow \Sigma A$ of
the coderivation $m$ correspond to maps $\tilde{m}_i: A^{\otimes
i}\rightarrow A$ of degree $i-2$. The map $\tilde{m}_1$ is the
original differential in $A$,  the map $\tilde{m}_2$ is a
multiplication up to homotopy and $\tilde{m}_i:A^{\otimes
i}\rightarrow A$ are  higher multiplications on $A$.

The space $Coder(T\Sigma A)\cong Hom(T\Sigma A,\Sigma A)$  is a
differential graded Lie algebra with respect to the (graded)
commutator. Let $m,n\in Hom(T\Sigma A,\Sigma A)$ so that $m
=(m_1,m_2,\ldots)$ and $n=(n_1,n_2,\ldots)$ where $m_i,n_i\in
Hom((\Sigma A^{\otimes i}),\Sigma A)$. Then the commutator of $m$
and $n$ is clearly determined by commutators of their components:
$[m_i,n_j]=m_i\circ n_j-(-1)^{|m_i||n_i|}n_i\circ m_i$.
Furthermore we have the following formula for the composition
$m_i\circ n_j\in Hom((\Sigma A)^{\otimes i+j-1},\Sigma A)$:
\begin{eqnarray}\label{big}(m_i\circ
n_j)[a_1|\ldots |a_{i+j-1}]\\=
\nonumber\sum_{k=0}^{i-1}(-1)^{|n_j| (|a_1|+\ldots+|a_k|+k)}
m_i[a_1|\ldots
|a_k|n_j[a_{k+1}|\ldots,|a_{k+j}]|a_{k+j+1}|\ldots|a_{i+j-1}]
\end{eqnarray}
\begin{defi} For two $A_\infty$-algebras $A$ and $C$ an
$A_\infty$-morphism (or $A_\infty$-map) $A\rightarrow C$ is a map
of differential graded coalgebras $f:{T}\Sigma A\rightarrow
{T}\Sigma C$ for which $f(T^0\Sigma A)=T^0\Sigma C=R$.
\end{defi}

It is clear that any $A_\infty$-map $T\Sigma A\rightarrow T\Sigma
C$ between two $A_\infty$-algebras $A$ and $C$ is a map of
filtered coalgebras. Furthermore a coalgebra map $f:T\Sigma
A\rightarrow \Sigma C$ is determined by the composition $T\Sigma
A\stackrel{f}{\rightarrow} T\Sigma C\rightarrow \Sigma C$ where
the second arrow is the canonical projection. Denoting the
components of the composite map by $f_i$ we see that $f$ is
determined by the collection $(f_1,f_2,\ldots)$ where $f_i:(\Sigma
A)^{\otimes i}\rightarrow \Sigma C$. The map $f$ could be
recovered from the collection $\{f_i\}$ as follows:
\begin{equation}\label{rec2}f[a_1|\ldots|a_n]=\sum[f_{i_1}[a_1|\ldots|a_{i_1}]|
\ldots|f_{i_k}[a_{i_{k-1}+1}|\ldots|a_{n}]]\end{equation} where
the summation is over all partitions $(i_1,\ldots,i_k)$ of $n$.

 The components $f_i$ of the $A_\infty$-map $f$
correspond to the maps $\tilde{f}_i:A^{\otimes i}\rightarrow C$ of
degrees $i-1$. The map $\tilde{f}_1:A\rightarrow C$ is a map of
complexes which is multiplicative up to higher homotopies provided
by $\tilde{f}_2,\tilde{f}_3,\ldots$.

We say that an $A_\infty$-map $f$ is a weak equivalence  if
$f_1:\Sigma A\rightarrow \Sigma C$ is a quasi-isomorphism of
complexes. Further we say that two $A_\infty$-algebras $A$ and $C$
are weakly equivalent if there is a chain of weak equivalences
$A\rightarrow A_1\leftarrow A_2\rightarrow\ldots\leftarrow A_n=C$.
\begin{rem}In fact it is possible to prove that for a weak
equivalence $A\rightarrow C$ there is always a weak equivalence
$C\rightarrow A$. This could be proved by constructing a closed
model category of cocomplete differential graded $R$-coalgebras
and identifying bar-constructions of $A_\infty$-algebras as
fibrant-cofibrant objects in this category. The discussion of such
matters would take us too far afield and we refer the reader to
\cite{lef} where this construction is carried out.\end{rem}

The weak equivalences $f=(f_1,f_2,\ldots)$ that we consider later
on in the paper will always have the property that the morphism
$f_1$ is invertible (this is always the case for so-called {\it
minimal} $A_\infty$-algebras, i.e. such that $m_1=0$). The
following proposition is a version of the formal implicit function
theorem.
\begin{prop} Let $f=(f_1,f_2,\ldots):TA\rightarrow TC$ be a map of
(filtered) coalgebras. Then $f$ is invertible if and only if
$f_1:A\rightarrow C$ is invertible.\end{prop}
\begin{proof}If $f$ is invertible with the inverse
$g=(g_1,g_2,\ldots)$ then clearly $g_1$ is the inverse to $f_1$.
Conversely suppose that $f_1$ is invertible. We will construct a
sequence of maps $g^n:TC\rightarrow TA$ such that $f\circ
g^n=id\mod O(n)$ as follows. Set $g^1=(f_1^{-1},0,0,\ldots)$.
Clearly $f\circ g^1=id\mod O(1)$. Now assume by induction that the
maps $g^n=(g^n_1,g^n_2,\ldots)$ have been constructed for $n\leq
k$. Then up to the terms of filtration $\geq k+1$ we have $f\circ
g^k=id +f_1\circ g^k_{k+1}+X$ where $X$ is some map $TA\rightarrow
TA$ having filtration $k$. Set
$g^{k+1}=g^k-g^k_{k+1}-f_1^{-1}\circ X$. Then $g^{k+1}$ agrees
with $g^k$ up to the terms of filtration $\geq k-1$ and has the
property that $f\circ g^{k+1}=id~\mod ~O(k+1)$. The sequence
$\{g^n\}$ clearly converges in the sense of the filtration on
$Hom(TC,TA$)and setting $g=lim_{n\rightarrow \infty}g^n$ we obtain
$f\circ g=id$. Similarly there exists a right inverse to $f$ so
$f$ is invertible.
\end{proof}
\section{Unital structures}
\begin{defi} An $A_\infty$ structure $m=(m_1,m_2,\ldots)$ on a cell complex
$A$ is called unital if there exists an element $1=1_A$ of degree
zero (called the unit of $A$) such that
$m_2[1_A|a]=a=(-1)^{|a|}m_2[a,1_A]$ for all $a\in A$ and such that
$m_i(a_1,\ldots,a_i)=0$ for all $i\neq 2$ if one of $a_i$ equals
$1_A$.  An $A_\infty$-morphism $f=(f_1,f_2,\ldots)$ between  two
unital $A_\infty$-algebras $A$ and $B$ is called unital if
$f_1[1_A]=[1_B]$ and $f_i[a_1|\ldots|a_i]=0$ for all $i\geq 2$ if
one of $a_i$ equals $1_A$.\end{defi} \begin{rem} Notice that
$m_2[1_A|a]=a=(-1)^{|a|}m_2[a|1_A]$ is equivalent to the more
customary $\tilde{m}_2(1_A,a)=\tilde{m}_2(a,1_A)=a$.\end{rem} From
now on we will use the term $A_\infty$-algebra for a unital
$A_\infty$-algebra and an $A_\infty$-morphism for a unital
$A_\infty$-morphism (unless indicated otherwise).

One important consequence of unitality (which will not be used in
this paper however) is that the complex $BA$ with the differential
determined by the collection $m=(m_1,m_2,\ldots)$ is exact for a
unital $A_\infty$-algebra $A$. We leave it to the interested
reader to check that the map $s:[a_1|a_2|\ldots |a_n]\rightarrow
[1|a_1|a_2|\ldots |a_n]$ is a contracting homotopy for $BA$.

The classification problem of $A_\infty$-algebras naturally leads
one to consider the group $Aut(T\Sigma A)$ of automorphisms of the
coalgebra $T\Sigma A$ where $A$ is a graded $R$-module. In the
unital case the relevant group is the group of {\it normalized}
automorphisms $\overline{Aut}(T\Sigma A)$ which we will now
define.
\begin{defi}Let $A$ be a free graded $R$-module with a
distinguished element $[1]\in \Sigma A$ of degree $1$. We call an
automorphism $f=(f_1,f_2,\ldots), f_i:(\Sigma A)^{\otimes
i}\rightarrow \Sigma A$ of the coalgebra $T\Sigma A$ normalized if
$f_1[1]=[1]$ while $f_i[a_1|\ldots|a_n]=0$ for $i>1$ if one of
$a_i$'s is equal to $1$. The set of normalized automorphisms will
be denoted by $\overline{Aut}(T\Sigma A)$.
\end{defi}
  The
set $\overline{Aut}(T\Sigma A)$ is in fact a group. Indeed using
the formula (\ref{rec2}) one sees immediately that the composition
of two normalized automorphisms is normalized. Therefore
$\overline{Aut}(T\Sigma A)$ is a subgroup in ${Aut}(T\Sigma A)$.

The concomitant notion to a normalized automorphism is that of a
normalized coderivation.
\begin{defi} A coderivation $\xi=(\xi_0,\xi_1,\ldots)\in Coder(T\Sigma A)$
will be called normalized if $\xi_i[a_1|\ldots|a_i]=0$ each time
one of $a_k=1$ for $ i=1,2,\ldots$. The set of all normalized
derivations is denoted by $\overline{Coder}(T\Sigma A)$.
\end{defi}
\begin{rem}Clearly the set $\overline{Coder}(T\Sigma A)$ forms a (graded)
Lie subalgebra in the Lie algebra $Coder(T\Sigma A)$. It is
natural to consider $\overline{Aut}(T\Sigma A)$ as the associated
Lie group. \end{rem} It is often extremely convenient to work in
the dual setting. Suppose that the element $1\in A$ can be
completed to a basis $\{1,y_i, i\in I\}$ of the $R$-module $A$.
The indexing set $I$ will be finite in our examples but need not
be in general.
\begin{rem}\label{loc}If our ground ring $R$ is local, than the above
assumption is always satisfied. Indeed let $\{e_i\}$ be a basis of
$A$ over $R$. Then $1=r_1e_1+\ldots+r_ne_n$. Clearly the element
$1$ remains nonzero after reducing modulo the maximal ideal in
$R$. Therefore one of the coefficients $r_1,\ldots,r_n$, say $r_1$
must be invertible in $R$. Then $1,e_2,e_3,\ldots$ form a basis in
$A$. Thus the assumption that $1$ can be completed to a basis in
$A$ is not really a restriction since we can always argue `one
prime at a time'.\end{rem}
 Then the $R$-module dual to the coalgebra $T\Sigma A$
(usually referred to as the cobar-construction)  is the algebra of
noncommutative power series in variables $\{\tau,{\bf t
}\}=\{\tau, t_1,t_2,\ldots\}$.  Here the elements $\tau,t_i$ form
the basis in $\Sigma A^*$ dual to $[1],[y_i]\in \Sigma A$:
$$(T\Sigma A)^\ast=k\langle\langle\tau, {\bf t}\rangle\rangle$$
Notice that $\tau$ has degree $-1$ whereas $|t_i|=-|y_i|-1$.

The algebra $R\langle\langle\tau, {\bf t}\rangle\rangle$ has a
linear topology where the fundamental system of neighborhoods of
$0$ is formed by those series whose constant term is $0$ and which
annihilate a finite dimensional submodule in $T\Sigma A$. It is
clear that $R\langle\langle\tau, {\bf t}\rangle\rangle$ is
Hausdorff and complete with respect to this topology.

Clearly the coalgebra endomorphisms of $T\Sigma A$ are in
one-to-one correspondence with continuous endomorphisms of the
algebra $R\langle\langle\tau, {\bf t}\rangle\rangle$ while
coderivations of $T\Sigma A$ are in one-to-one correspondence with
continuous derivations of $R\langle\langle\tau, {\bf
t}\rangle\rangle$. A continuous endomorphism $f$  of
$R\langle\langle\tau, {\bf t}\rangle\rangle$ is specified by its
values on $\tau$, which is a series $G(\tau, \bf t)$ of degree
$-1$ and on $t_i$'s which are series $F_i(\tau,\bf t)$ whose
degree equals that of $t_i$.  So $f$ corresponds to a collection
of power series of the form $(G(\tau, {\bf t}),F_1(\tau, {\bf
t}),F_2(\tau, {\bf t})\ldots)$. (Observe that if the indexing set
$I$ is infinite then continuity imposes certain restrictions on
the collection $G({\bf t}), F_1({\bf t}), F_2({\bf t}),\ldots)$).
The composition of endomorphisms corresponds to substitution of
power series. Similarly any continuous derivation $\xi$ could be
uniquely represented in the form $\xi=A(\tau,{\bf
t})\partial_\tau+ \sum_{i}B_i(\tau,{\bf t})\partial_{t_i}$. Here
$\partial_\tau$ and $\partial_{t_i}$ are standard derivations
corresponding to the coordinates $\tau,t_i$.
\begin{defi}A continuous derivation of
$R\langle\langle\tau,{\bf t}\rangle\rangle$ is called normalized
if the corresponding coderivation of $T\Sigma A$ is normalized. We
will denote the set of normalized derivations of
$R\langle\langle\tau,{\bf t}\rangle\rangle$ by
$\overline{Der}(R\langle\langle\tau,{\bf t}\rangle\rangle)$.
Similarly we call a continuous automorphism of
$R\langle\langle\tau,{\bf t}\rangle\rangle$ normalized if such is
the corresponding automorphism of $T\Sigma A$. The set of
normalized automorphisms of $R\langle\langle\tau,{\bf
t}\rangle\rangle$ will be denoted by
$\overline{Aut}(R\langle\langle\tau,{\bf t}\rangle\rangle)$.
\end{defi}

Recall that the space $Coder(T\Sigma A)$ has a filtration
$O(-1)\supset O(0)\supset\ldots$ where $O(n)$ consists of those
coderivations $\xi=(\xi_0,\xi_1,\ldots)$ for which
$\xi_1=\xi_2=\ldots=\xi_n=0$. Then the space of (continuous)
derivations of $(T\Sigma A)^\ast$ acquires filtration so that the
derivation $A(\tau,{\bf t})\partial_{\tau}+\sum_{i\in
I}B_i(\tau,{\bf t})\partial_{t_i}$ has weight $\geq n$ if and only
if the expressions $A(\tau,{\bf t}),B_1(\tau,{\bf
t}),B_2(\tau,{\bf t})\ldots)$ do not contain terms of degree $\leq
n$. We will still denote the collection of elements of weight
$\geq n$ by $O(n)$.
\begin{prop}\label{doi}$(i)$ Any normalized derivation $\xi$ of
$R\langle\langle\tau,{\bf t}\rangle\rangle$ has the form
$$\xi=A({\bf t})\partial_{\tau}+\sum_{i\in I}B_i({\bf
t})\partial_{t_i}.$$ $(ii)$ Any unital $A_\infty$ structure $m$ on
$A$ corresponds to a derivation $m^\ast$ of $(T\Sigma A)^\ast$ of
the form $$m^\ast=(A({\bf t})+\tau^2)\partial_{\tau}+\sum_{i\in
I}([\tau,t_i]+B_i({\bf t}))\partial_{t_i}.$$ where the series
$A({\bf t}), B({\bf t})$ have vanishing constant terms.
\end{prop}
\begin{proof}Denote by $\langle,\rangle$ the $R$-linear pairing between
$(T\Sigma A)$ and $(T\Sigma A)^\ast$. Associated to a homogeneous
endomorphism $T$ of the $R$-module $(T\Sigma A)$ is the
endomorphism $T^\ast$ of $(T\Sigma A)^\ast$ for which
\begin{equation}\label{dua}\langle T(a),b\rangle=(-1)^{|a||T|}\langle a,
T^\ast(b)\rangle\end{equation}
 The rest is just a routine exercise in dualization using (\ref{dua}) which
we can safely leave to the reader. Note that the quadratic term
$\sum_{i\in I}[\tau,t_i]\partial_{t_i}+\tau^2\partial_{\tau}$
corresponds to the identities
$m_2[1|y_i]=(-1)^{|y_i|}m_2[y_i|1]=[y_i]$ and $m[1|1]=[1]$. The
condition that $A({\bf t}), B({\bf t})$ have vanishing constant
terms means that $m^\ast\in O(0)$.
\end{proof}
\begin{rem}Of course not every derivation $\xi$ of $(T\Sigma A)^\ast$ of
the form $\xi=(A({\bf t})+\tau^2)\partial_{\tau}+\sum_{i\in
I}([\tau,t_i]+B_i({\bf t}))\partial_{t_i}$ is an $A_\infty$
structure. The condition which specifies an $A_\infty$ structure
is $\xi\circ \xi=0$ (or, equivalently, $[\xi,\xi]=0$ if $R$ has no
$2$-torsion). Also the condition that $m^\ast$ has degree $-1$
puts further restrictions on $A({\bf t})$ and $B_i({\bf t})$. For
example if the variables $t_i$ have even degrees then all $B_i$'s
necessarily vanish.
\end{rem}
\begin{rem}It is easy to check  that the derivation $\sum_{i\in
I}([\tau,t_i])\partial_{t_i}+\tau^2\partial_{\tau}$ can be
compactly written as $ad\tau-\tau^2\partial_{\tau}$ where
$ad\tau(?):=[\tau,?]$. Therefore the formula for $m^\ast$ could be
written as $$m^\ast=A({\bf t})\partial_{\tau}+\sum_{i\in
I}B_i({\bf t})\partial_{t_i}+ad\tau-\tau^2\partial_\tau.$$
\end{rem}
Similarly we could translate the notion of a normalized
automorphism to the dual setting. Consider the continuous
endomorphism of $R\langle\langle\tau,{\bf t}\rangle\rangle$
corresponding to the collection $(G,{\bf F}):=(\tau+G({\bf
t}),F_1({\bf t}),F_2({\bf t}),\ldots)$ of power series without
constant terms. Here we require that ${\bf F}({\bf t} )=(F_1({\bf
t}),F_2({\bf t}),\ldots):R\langle\langle{\bf
t}\rangle\rangle\rightarrow R\langle\langle{\bf t}\rangle\rangle$
be invertible with inverse ${\bf F}^{-1}({\bf t})$. Then clearly
$(G,{\bf F} )$ is invertible and $(G,{\bf F} )^{-1}=(-G({\bf
F}^{-1}),{\bf F }^{-1})$. Moreover such endomorphisms form a
subgroup of all continuous automorphisms of
$R\langle\langle\tau,{\bf t}\rangle\rangle$. Then we have the
following result whose proof is similar to part $(i)$ of
Proposition \ref{doi}.
\begin{prop} The group of continuous automorphisms of
$R\langle\langle\tau,{\bf t}\rangle\rangle$ consisting of pairs
$(G,{\bf F})$ as above is isomorphic to
$\overline{Aut}(R\langle\langle{\tau,\bf t}\rangle\rangle)$.
\end{prop}
\begin{rem}The condition that a multiplicative automorphism
necessarily has degree zero puts certain restrictions on ${\bf F}$
and $G$. For example if all variables $t_i$ have even degrees then
$G({\bf t})=0$.\end{rem} \begin{rem}It is illuminating to consider
the unit map $R\rightarrow A$ from the point of view of the cobar
construction. Observe that the canonical structure of an
associative algebra on $R$ corresponds to the derivation
$\tau^2\partial_\tau$ of the power series ring $R[[\tau]]$. Then
the unit map $R\rightarrow A$ considered as an $A_\infty$-map is
the map of cobar constructions $$(T\Sigma A)^\ast
=R\langle\langle\tau,{\bf t
}\rangle\rangle\stackrel{i}{\longrightarrow} (T\Sigma R)^\ast=
R[[\tau]]$$ where $i(\tau)=\tau$ and $i({\bf t})= 0$. The
unitality condition ensures that $i$ is a map of dga's. Further
the maps of dga's $R[[\tau]]\rightarrow R\langle\langle\tau, {\bf
t}\rangle\rangle$ should be considered as  `$A_\infty$-points' of
$A$. The existence of $A_\infty$-points is a subtle question in
general and we hope to return to it in in the future. If the
$A_\infty$ structure $m^\ast$ has the form
$m^\ast=ad\tau-\tau^2\partial_\tau+\sum_{i\in I}B({\bf
t})\partial_{t_i}$ then the map $\epsilon:R[[\tau]]\rightarrow
R\langle\langle\tau, {\bf t}\rangle\rangle:\epsilon(\tau)=\tau$ is
a `canonical' $A_\infty$-point of $A$.
\end{rem} Next observe that the group $Aut(T\Sigma A)$ acts on
the set of coderivations of $\Sigma TA$ according to the formula
$f:m\rightarrow m^f=f\circ m\circ f^{-1}$ for $m\in Coder(\Sigma
TA)$ and $f\in Aut(\Sigma TA)$. Obviously if $m\circ m=0$ then
$m^f\circ m^f=0$ so $Aut(\Sigma TA))$ acts on the set of
(nonunital) $A_\infty$ structures on $A$. It turns out the the
group $\overline{Aut}(T\Sigma
A)=\overline{Aut}(R\langle\langle{\tau,\bf t}\rangle\rangle)$ acts
on the set of {\it unital} $A_\infty$ structures.

Denote by $(A,{\bf B})$ the derivation of
$R\langle\langle\tau,{\bf t}\rangle\rangle$ corresponding to a
unital $A_\infty$ structure: $$(A,{\bf B})=(A({\bf
t})+\tau^2)\partial_{\tau}+\sum_{i\in I}([\tau,t_i]+B_i({\bf
t}))\partial_{t_i}.$$ \begin{prop}\label{act}The group
$\overline{Aut}(R\langle\langle{\tau,\bf t}\rangle\rangle)$ acts
on the right on the set of unital $A_\infty$ structures according
to the formula
\begin{eqnarray}\label{we}(A,{\bf B})\ast(G, {\bf F})=(G,
{\bf F})\circ (A,{\bf B})\circ (G, {\bf F})^{-1}\\ \nonumber
=(A({\bf F(t)})-G({\bf t})^2+\sum_{j\in I}[B_j({\bf F}({\bf
t}))\partial_{t_j}G({\bf F^{-1})]}({\bf F(t)}),\sum_{i,j\in
I}([G({\bf t}),t_i]+(B_j({\bf F}( {\bf t}))\partial_{t_j}{\bf
F}^{-1})({\bf F}(t_i)).\end{eqnarray}
\end{prop}
\begin{proof}This is one of the examples where the use of the dual language
leads to considerable simplifications; the relatively painless
calculations below become exceedingly gruesome when performed in
terms of coderivations of the coalgebra $T\Sigma A$. We compute:
\begin{eqnarray}\label{efr}
\nonumber((A,{\bf B})\ast (G, {\bf F}))(t_i)=(G, {\bf F})\circ
(A,{\bf B})\circ (G, {\bf F})^{-1}(t_i)=(G, {\bf F})\circ (A,{\bf
B})({ \bf F}^{-1}( t_i))\\ \nonumber=(G, {\bf F})([\tau,{\bf
F}^{-1}( t_i)]+\sum_{j\in I}B_j{(\bf t})\partial_{{ t}_j}{\bf
F}^{-1}( t_i))\\  =[\tau,t_i]+[G({\bf t}),t_i]+\sum_{j\in
I}(B_j({\bf F}( {\bf t}))\partial_{t_j}{\bf F}^{-1})({\bf
F}(t_i))\end{eqnarray} Further
\begin{eqnarray}\nonumber((A,{\bf B})\ast (G, {\bf
F}))(\tau)=(G, {\bf F})\circ (A,{\bf B})\circ (G, {\bf
F})^{-1}(\tau)=(G, {\bf F})\circ (A,{\bf B})(\tau -G({\bf
F}^{-1}({\bf t})))\\ \nonumber=(G, {\bf F})(\tau^2+A({\bf
t})-[\tau, G({\bf F}^{-1}({\bf t}))]-\sum_{j\in I}B_j({\bf
t})\partial_{t_j}G({\bf F}^{-1}({\bf t})))\\
\nonumber=(\tau+G({\bf t}))^2+A({\bf F}({\bf t}))-[\tau+G({\bf
t}),G({\bf t})]+\sum_{j\in I}[B_j({\bf F}({\bf
t}))\partial_{t_j}G({\bf F^{-1})]}({\bf F(t)}).\end{eqnarray}
Since $G({\bf t})$ and $\tau$ have odd degrees we have the
equalities $(\tau+G({\bf t}))^2=\tau^2+G({\bf t })^2+[\tau,G({\bf
t})]$ and $[G({\bf t}),G({\bf t})]=2G({\bf t})^2$.   It follows
that
\begin{equation}\label{ef}((A,{\bf B})\ast (G, {\bf F}))(\tau)=\tau^2+
A({\bf F}({\bf t}))-G({\bf t})^2+\sum_{j\in I}[B_j({\bf F}({\bf
t}))\partial_{t_j}G({\bf F^{-1})]}({\bf F(t)})\end{equation}The
formula (\ref{we}) is a consequence of (\ref{efr}) and (\ref{ef})
and our proposition is proved.
\end{proof}
\begin{rem}The above proposition has two parts:
the statement that the group
$\overline{Aut}(R\langle\langle{\tau,\bf t}\rangle\rangle)$ acts
on the set of unital $A_\infty$ structures and an explicit formula
for this action. While the formula clearly requires the assumption
that $1$ can be completed to an $R$-basis in $A$, the statement
about the group action is valid without this assumption. The proof
of this could be deduced from Remark \ref{loc} using standard
localization techniques.\end{rem}
\begin{rem}Proposition \ref{act} admits the following infinitesimal analogue:
if $m$ is a unital $A_\infty$ structure and $\xi$ is a normalized
coderivation of $T\Sigma A$ then the commutator $[\xi,m]$ is also
normalized. This can be interpreted as saying that the normalized
Hochschild  cochains of a unital $A_\infty$-algebra form a
subcomplex with respect to the Hochschild differential, cf. next
section of the present paper.\end{rem}
\begin{rem}\label{act1}We have seen that if the variables $t_i$ all have even degrees
then $B_i=0$ and $G=0$. In other words the group of normalized
automorphisms is just the group of formal power series ${\bf
F}({\bf t})$ under composition and a unital $A_\infty$ structure
corresponds to the derivation of the form $A({\bf
t})\partial_{\tau}$. The formula (\ref{we}) in this case takes an
especially simple form: $A\ast{\bf F}=A({\bf F})$.
\end{rem}

The next result we are going to discuss requires a certain
knowledge of {\it operads}. We do not intend to discuss this
subject in detail here and refer the interested reader to the nice
exposition in \cite{Vor}. An $A_\infty$-algebra is in fact an
algebra over a certain operad $\mathcal{A_\infty}$ in the category
of differential graded $R$-modules, sometimes called the
\emph{Stasheff operad}. The operad $\mathcal{A}_\infty$ maps into
another operad $\mathcal{A}ss$ whose algebras are strictly
associative differential graded $R$-algebras and this map is a
quasi-isomorphism. In particular any differential graded algebra
is an $A_\infty$-algebra.
\begin{prop}\label{dif}There is a functor that assigns to each unital
 $A_\infty$-algebra $A$ a strictly associative differential graded algebra
 $\tilde{A}$ which is weakly equivalent to $A$.\end{prop}
\begin{proof} We will only give a sketch following \cite{KM}, $V.1.7$.
 Associated to any operad is a monad
having the same algebras. Denote the monad in the category of
complexes of $R$-modules associated to $\mathcal{A}ss$ by $C$ and
the one associated to $\mathcal{A_\infty}$ by $C_\infty$. Then
there is a canonical map of monads $C\rightarrow C_\infty$.
Consider the following maps of $\mathcal{A_\infty}$-algebras
\begin{equation}\label{ok}A\leftarrow
B({C_\infty},{C_\infty},A)\rightarrow
B(C,{C_\infty},A)\end{equation} Here $B(-,-,A)$ stands for a
two-sided monadic bar construction. Both maps in (\ref{ok}) are
homology isomorphisms and our proposition is proved.\end{proof}

\section{Hochschild Cohomology of $A_\infty$-algebras.}

Let $A$ be an $A_\infty$-algebra. Consider the graded Lie algebra
$Coder(BA)$ of all  coderivations of the coalgebra $BA=T\Sigma A$.
There is a preferred coderivation $m:BA\rightarrow BA$ of degree
$-1$ which is given by the $A_\infty$ structure on $A$. We will
define a differential $\partial$ on $Coder(BA)$ by the formula
$\partial(f)=[f,m]$ where the right hand side is the (graded)
commutator of two coderivations $f$ and $m$. The condition $m\circ
m=0$ implies  that $\partial\circ\partial=0$.
\begin{defi}The complex $C^\ast(A,A):=Coder(BA)$ with the differential $\partial$
is called the Hochschild complex  of an $A_\infty$-algebra $A$.
Its cohomology $H^\ast(A,A)$ is called the Hochschild cohomology
of $A$ with coefficients in itself.\end{defi}
\begin{rem}Since the coderivation $m$ has weight $\geq 0$
the differential on $C^\ast(A,A)$ agrees with the filtration on
$BA$ in the sense that $d(O(n))\subset O(n)$.\end{rem} Recall that
since $BA$ is cofree in the category of cocomplete coalgebras
there is a natural identification $C^\ast(A,A)\cong Hom(BA,\Sigma
A)$ which we will use without explicitly mentioning. Using the the
formula (\ref{rec1}) one can recover the coderivation of
$BA=T\Sigma A$ from its components $c_k\in Hom((T\Sigma
A)^{\otimes k },\Sigma A)\subset Hom(BA,\Sigma A)$.

 We will now introduce the
normalized Hochschild complex for $A_\infty$-algebras which is
smaller and easier to compute with.
\begin{defi}Let $A$ be an $A_\infty$-algebra.
Then a Hochschild cochain $c\in Hom(T\Sigma A,\Sigma A)$ is called
normalized if $c$ is  normalized  as a coderivation of $BA$.
\end{defi} It is easy to check using (\ref{big}) that
the normalized cochains form a subcomplex of the Hochschild
complex. We will denote this subcomplex by $\bar{C}^\ast(A,A)$.
\begin{theorem}\label{norm} Let $A$ be an $A_\infty$-algebra. Then
there is a chain deformation retraction of $C^\ast(A,A)$ onto the
subcomplex $\bar{C}^\ast(A,A)$. In particular both complexes have
the same cohomology.\end{theorem}
\begin{proof}
The proof is similar to that of the classical theorem of
Eilenberg-MacLane on normalized simplicial modules. Note that this
theorem cannot be applied directly since the Hochschild cohomology
of an $A_\infty$-algebra \emph{is not} a cohomology of a
simplicial object.  The resulting calculations in the $A_\infty$
context are considerably more involved.

Let us call a cochain $c\in C^\ast(A,A)$ $i$-normalized if $c$
vanishes each time one if its first $i$ arguments is equal to $1$.
Then $c$ is normalized if and only if it is $i$-normalized for all
$i$.

We define a sequence of cochain maps $h_i:C^\ast(A,A)\rightarrow
C^{\ast}(A,A)$ as follows. Let $c\in Hom((\Sigma A)^{\otimes
n},\Sigma A)$ for some $n$ and consider the cochain $s_i(c)\in
Hom((\Sigma A)^{\otimes n-1},\Sigma A)$ defined by the formula
$$s_i(c)[a_1|\ldots|a_{n-1}]=(-1)^{|a_1|+\ldots+|a_{l}|+i+1}
c[a_1|\ldots|a_{i}|1|a_{i+1}|\ldots|a_{n-1}].$$ Extending by
linearity we define $s_i$ on the whole $C^{\ast}(A,A)$. Then set
$h_i(c):=c-\partial(s_i(c))-s_i(\partial c)$. We claim that $h_i$
takes an $i$-normalized Hochschild cochain to an $i+1$-normalized
cochain. Indeed, let $c$ be an $i$-normalized cochain. We could
assume without loss of generality that $c\in Hom((\Sigma
A^{\otimes n}), \Sigma A)$ for some $n>i$. We want to show that
\begin{equation}\label{er}c(?)=[s_i(c),m_k](?)+s_i[c,m_k](?)\end{equation}
 for any $k$ as long as
the $i+1$st argument in $?$ is $1$.

Notice that the left hand side of (\ref{er}) is only nonzero if
the number of arguments in $?$ is $n$ whereas the right hand side
of (\ref{er}) is nonzero if the number of arguments is $n-k+2$.
Therefore we need to consider the cases $k= 2$ and $k\neq 2$
separately. For $k= 2$ we have
\begin{eqnarray}\nonumber (m_2\circ s_ic)[a_1|\ldots|a_{n}]
=(-1)^{|a_1|+\ldots+|a_i|+i+1}
m_2[c[a_1|\ldots|a_i|1|a_{i+1}|\ldots a_{n-1}]a_n]\\ \nonumber \pm
m_2[a_1|c[a_2|\ldots|a_{i+1}|1|a_{i+2}|\ldots|a_n].
\end{eqnarray}
 Setting $a_{i+1}=1$ and taking into
account that $c$ is $i$-normalized we obtain $$(m_2\circ
s_ic)[a_1|\ldots|a_i|1|a_{i+2}| \ldots|a_{n}]
=(-1)^{|a_1|+\ldots+|a_i|+i+1}m_2
[c[a_1|\ldots|a_i|1|1|a_{i+2}|\ldots a_{n-1}]a_n].$$

Similarly we obtain $$s_i(m_2\circ c
)[a_1|\ldots|a_i|1|a_{i+2}|\ldots|a_n]
=(-1)^{|a_1|+\ldots+|a_i|+i+1}m_2
[c[a_1|\ldots|a_i|1|1|a_{i+2}|\ldots a_{n-1}]a_n].$$ It follows
that
\begin{equation}\label{pou}(m_2\circ
s_ic-s_i(m_2\circ c
))[a_1|\ldots|a_i|1|a_{i+2}|\ldots|a_n]=0\end{equation}
 Taking into account the identities $m_2[a_i|1]=(-1)^{|a_i|}[a_i]$ and
$m_2[1|a_{i+1}]=[a_{i+1}]$ we have
\begin{eqnarray}\label{rt2}\nonumber (s_ic\circ m_2)
[a_1|\ldots|a_n]=\sum_{l=0}^i\pm
s_ic[a_1|\ldots|a_l|m_2[a_{l+1}|a_{l+2}]|a_{l+3}|\ldots|\ldots|a_n]\\
\nonumber
+\sum_{l=i+1}^{n-2}(-1)^{|a_1|+\ldots+|a_{l}|+l}s_ic[a_1|\ldots|a_i|a_{i+1}|
\ldots|a_l|m_2[a_{l+1}|a_{l+2}]|a_{l+3}|\ldots|a_n].
\end{eqnarray}
 After substituting $a_{i+1}=1$ the term having the sign $\pm$ in
front of it vanishes and we get
\begin{eqnarray}\nonumber (s_ic\circ
m_2)[a_1|\ldots|a_i|1|a_{i+1}|\ldots|a_n] \\=\nonumber
\sum_{l=i+1}^{n-2}(-1)^{|a_1|+\ldots+|a_{l}|+l+|a_1|+\ldots+|a_i|+i+1}
c[a_1|\ldots|a_i|1|1|\ldots|a_l|m_2[a_{l+1}|a_{l+2}]|a_{l+3}|\ldots|a_n].
\end{eqnarray}
 And similarly
 \begin{eqnarray}\nonumber
s_i(c\circ
m_2)[a_1|\ldots|a_i|1|\ldots|a_n]=c[a_1|\ldots|a_i|1|a_{i+1}|\ldots|a_n]\\
\nonumber+\sum_{l=i+1}^{n-2}(-1)^{|a_1|+\ldots+|a_{l}|+l+1
+|a_1|+\ldots+|a_i|+i}
c[a_1|\ldots|a_i|1|1|a_{i+2}|\ldots|a_l|m_2[a_{l+1}|a_{l+2}]|a_{l+3}|\ldots|a_n]
\end{eqnarray}
Therefore \begin{equation}\label{tr1}(s_ic\circ m_2+s_i(c\circ
m_2))[a_1|\ldots|a_i|1|a_{i+1}|\ldots|a_n]=c[a_1|\ldots|a_i|1|a_{i+1}|\ldots|a_n]
\end{equation}
 Taking into account that $|s_ic|=|c|+1$ we obtain from (\ref{tr1}), and (\ref{pou})
\begin{eqnarray}\nonumber
([s_ic,m_2]+s_i[c,m_2])[a_1|\ldots|a_i|1|a_{i+1}|\ldots|a_n]\\
\nonumber =( s_ic\circ m_2-(-1)^{|c|+1}m_2\circ s_ic+s_i(c\circ
m_2)- (-1)^{|c|}s_i(m_2\circ
c))[a_1|\ldots|a_i|1|a_{i+1}|\ldots|a_n]\\ \nonumber
=c[a_1|\ldots|a_i|1|a_{i+1}|\ldots|a_n].
\end{eqnarray}

Now let $k\neq 2$. We have: \begin{eqnarray}\nonumber(m_k\circ
s_ic)[a_1|\ldots|a_i|1|a_{i+2}|a_{n+k-2}]=
m_k[s_ic[a_1|\ldots|a_{n-1}]a_n|\ldots|a_{n+k-2}]\\ \nonumber
=(-1)^{|a_1|+\ldots+|a_{i}|+i+1}
m_k[c[a_1|\ldots|a_i|1|1|a_{i+2}|\ldots|a_{n-1}]|a_n|\ldots|a_{n+k-2}]
\end{eqnarray}
(the remaining terms in the expansion for $(m_k\circ
s_ic)[a_1|\ldots|a_i|1|a_{i+2}|a_{n+k-2}]$ vanish because the
Hochschild cochains $s_ic$ is $i$-normalized and $m_k$ is
normalized). Likewise
\begin{eqnarray}\label{re1}\nonumber
s_i(m_k\circ c)[a_1|\ldots|a_i|1|a_{i+1}|\ldots|a_{n+k-2}]\\
\nonumber= (-1)^{|a_1|+\ldots+|a_i|+i+1}m_k\circ
c[a_1|\ldots|a_i|1|1|a_{i+2}|\ldots|a_{n+k-2}]\\ \nonumber
=(-1)^{|a_1|+\ldots+|a_{i}|+i+1}
m_k[c[a_1|\ldots|a_i|1|1|a_{i+2}|\ldots|a_{n-1}]|a_n|\ldots|a_{n+k-2}]
\end{eqnarray}
It follows that
\begin{equation}\label{pou1}(m_k\circ
s_ic-s_i(m_k\circ c
))[a_1|\ldots|a_i|1|a_{i+2}|\ldots|a_n]=0.\end{equation} Further
\begin{eqnarray}\nonumber (s_ic\circ
m_k)[a_1|\ldots|a_i|1|a_{i+1}|\ldots|a_{n+k-2}]=\\ \nonumber
\sum_{l=i}^{n-2}(-1)^{|a_1|+\ldots+|a_l|+l}
s_ic[a_1|\ldots|a_i|1|a_{i+1}|\ldots
|a_l|m_k[a_{l+1}|\ldots|a_{l+k}]|a_{l+k+1}|\ldots|a_{k+n-2}]\\
\nonumber =\sum_{l=i}^{n-2}\epsilon_l
c[a_1|\ldots|a_i|1|a_{i+1}|\ldots
|a_l|m_k[a_{l+1}|\ldots|a_{l+k}]|a_{l+k+1}|\ldots|a_{k+n-2}]
\end{eqnarray} where $\epsilon_l=(-1)^{|a_1|+\ldots+|a_i|+i+1+|a_1|+\ldots+|a_l|+l}$
and similarly
\begin{eqnarray}\nonumber s_i(c\circ
m_k)[a_1|\ldots|a_i|1|a_{i+1}|\ldots|a_{n+k-2}]=
(-1)^{|a_1|+\ldots+|a_i|+i+1}c\circ
m_k[a_1|\ldots|a_i|1|a_{i+1}|\ldots|a_{k+n-2}]\\ \nonumber=
\sum_{l=i+1}^{n-2}(-\epsilon_l)
c[a_1|\ldots|a_i|1|1|a_{i+2}|\ldots
|a_l|m_k[a_{l+1}|\ldots|a_{l+k}]|a_{l+k+1}|\ldots|a_{k+n-2}].
\end{eqnarray}
Therefore \begin{equation}\label{pou2}(s_ic\circ m_k+s_i(c\circ
m_k))[a_1|\ldots|a_i|1|a_{i+1}|\ldots|a_{n+k-2}]=0\end{equation}
Finally (\ref{pou1}) and (\ref{pou2}) imply that $[s_ic,
m_k]+s_i[c ,m_k]=0.$

This proves (\ref{er}) and, therefore, our claim that $h_i$ takes
$i$-normalized cochains into $i+1$-normalized cochains.  It
follows that the composition $\ldots\circ h_l\circ h_{l-1}\ldots
h_0$ takes an arbitrary cochain $c\in C^\ast(A,A)$ into a
normalized cochain and exhibits the subcomplex $\bar{C}(A,A)$ of
normalized cochains as a chain deformation retract of
$C^\ast(A,A)$.
\end{proof}
\begin{rem}Part of Theorem \ref{norm} could be interpreted as
saying that if the cochain $c\in C^\ast(A,A)$ has the property
that $[c,m]$ belongs to the Lie subalgebra of normalized cochains
then there exists a normalized cochain $c^\prime$ for which
$[c,m]=[c^\prime,m]$. This result has the following globalization
proved in \cite{lef}: if two minimal $A_\infty$ structures are
equivalent through a nonunital $A_\infty$-morphism, then they are
equivalent also through a unital one. It would be interesting to
deduce this result from Theorem \ref{norm} (the proof in the cited
reference uses obstruction theory).
\end{rem}
\begin{rem} Let us call an $A_\infty$-algebra $A$ \emph{homotopy
unital}
if there exists an element $1\in A$ of degree $0$ for which
$m_1[1]=0$, and $m_2[1|a]=(-1)^{|a|}m_2[a|1]$ for any $a\in A$. It
is easy to see that (in contrast with strict unitality) weak
equivalences preserve homotopy unitality for minimal
$A_\infty$-algebras. Then it is proved in \cite{lef} that any
minimal homotopy unital $A_\infty$-algebra is weakly equivalent to
a (strictly) unital one. This result combined with the previous
remark and Proposition \ref{act} shows that the classification
problem of (minimal) homotopy unital $A_\infty$-algebras up to a
nonunital weak equivalence is equivalent to classification of
unital $A_\infty$-algebras up to a unital weak equivalence.
\end{rem}
For the next result we will need a slightly more general
definition of Hochschild cohomology than the one already given.
Let $(A,m_A)$, $(C,m_C)$ be two $A_\infty$-algebras and
$i:BA\rightarrow BC$ an $A_\infty$-morphism between them. We say
that a map $f:BA\rightarrow BC$ is a coderivation of the coalgebra
$BA$ with values in the coalgebra $BC$ if the following diagram is
commutative: $$\xymatrix{BA\ar[d]^f\ar[r]^-{\Delta_{BA}}&BA\otimes
BA\ar[d]^{f\otimes i+i\otimes f}\\
BC\ar[r]^-{\Delta_{BC}}&BC\otimes BC}$$ Here $\Delta_{BA}$ and
$\Delta_{BC}$ denote the diagonals in the coalgebras $BA$ and
$BC$. Then the space $Coder(BA,BC)$ becomes a complex with the
differential $df=m_C\circ f-(-1)^{|f|}f\circ m_A$. We will denote
this complex by $C^\ast(A,C)$.

Now let $c\in C^\ast(A,A)$ be a Hochschild cochain. Define the
cochain $i_\ast(c)\in C^\ast(A,C)$ by the formula
$i_\ast(c)=i\circ c:BA\rightarrow BC$. Likewise for a cochain
$c^\prime\in C^\ast(C,C)$ define the cochain $i^\ast(c^\prime)\in
C^\ast(A,C)$ by the formula $i^\ast(c^\prime)=c^\prime\circ i$. It
is straightforward to check $i_\ast$ and $i^\ast$ give maps of
cochain complexes: $$i_\ast:C^\ast(A,A)\rightarrow
C^\ast(A,C)\leftarrow C^\ast(C,C):i^\ast.$$
\begin{prop}\label{lop}For two weakly equivalent $A_\infty$-algebras $A$ and
$C$ their Hochschild complexes $C^\ast(A,A)$ and $C^\ast(C,C)$ are
quasi-isomorphic as complexes of $R$-modules. In particular,
$H^\ast(A,A)\cong H^\ast(C,C)$.\end{prop}
\begin{proof} Let $i:BA\rightarrow BC$ be an $A_\infty$-morphism
establishing a weak equivalence between $A$ and $C$. Since the
cochain map $i_\ast:C^\ast(A,A)\rightarrow C^\ast(A,C)$ is a
filtered map it induces a map on associated spectral sequences.
Since $i$ induces a quasi-isomorphism $A\rightarrow C$ we see that
$i_\ast$ induces an isomorphism of the $E_1$-terms of the
corresponding spectral sequences. Therefore $i_\ast$ is itself a
quasi-isomorphism. Similar considerations show that the cochain
map $i^\ast:C^\ast(C,C)\rightarrow C^\ast(A,A)$ is a
quasi-isomorphism and our proposition is proved.
\end{proof}
\begin{rem} In general the complex $C^\ast(A,A)$ as
well as its cohomology $H^\ast(A,A)$ is not functorial with
respect to $A$. It is possible to define the Hochschild complex
$C^*(A,M)$ of an $A_\infty$-algebra with coefficients in a
$A_\infty$-bimodule $M$, cf. for example, \cite{GJ1}. Then
$C^*(A,M)$ is contravariant in the variable $A$ and covariant in
the variable $M$. However we don't need such level of generality
here and the discussion of $A_\infty$-bimodules would take us too
far afield.
\end{rem} Now let $A$ be an $A_\infty$-algebra.
Propositions \ref{lop} and \ref{dif} shows that the complex
$C^\ast(A,A)$ is quasi-isomorphic to the complex
$C^\ast(\tilde{A},\tilde{A})$ where $\tilde{A}$ is a differential
graded (unital) algebra weakly equivalent to $A$. The complex
$C^\ast(\tilde{A},\tilde{A})$ is the usual Hochschild complex of
the  dga  $\tilde{A}$ and it is well-known that it possesses
itself a structure of a homotopy commutative dga; something that
we did not see from the point of view of the $A_\infty$-algebra
$A$. The Hochschild complex of a dga admits a different (but of
course equivalent) description. Namely, we can define the complex
$C^\ast(\tilde{A},\tilde{A})$ as an object in the derived category
of $\tilde{A}\otimes\tilde{A}^{op}$-modules:
$$C^\ast(\tilde{A},\tilde{A}):=
RHom_{\tilde{A}\otimes^L\tilde{A}^{op}}(\tilde{A},\tilde{A}).$$
Here $\tilde{A}^{op}$ is the differential graded algebra having
the same underlying $R$-module and differential as $\tilde{A}$ but
the opposite multiplication. Since $\tilde{A}$ has the same
homology algebra as $A$ we get the following result:
\begin{prop}\label{ss} There exists a spectral sequence of $R$-modules
 $$Ext^{\ast\ast}_{H_\ast(A_\ast\otimes^L
A_\ast^{op})}(H_\ast(A),H_\ast(A))\Longrightarrow H^\ast(A,A).$$
It is of standard cohomological type, lies in the right half plane
and converges conditionally.\end{prop}
\section{Moore algebras}
In this section we introduce and study a class of
$A_\infty$-algebras which will be called Moore $A_\infty$-algebras
or just Moore algebras. The terminology comes from stable homotopy
theory - a Moore algebra is analogous to the Moore spectrum which
is a cofibre of the map $S\stackrel{p}{\rightarrow}S$ where $S$ is
the sphere spectrum. In some sense Moore algebras are the simplest
nontrivial examples of $A_\infty$-algebras which are not
differential graded algebras.

\begin{defi}An $A_\infty$-algebra over a commutative evenly graded
ring $R$ is called a Moore algebra if its underlying complex is
$A=\{\Sigma^dR\stackrel{\partial}{\rightarrow}R\}$ for some
differential $\partial$ . The integer $d$ is called the degree of
$A$.\end{defi} Obviously the generator in degree $0$ is $1\in R$.
We will denote the generator in $\Sigma^dR$ by $y$, so $|y|=d+1$.
The structure of an $A_\infty$-algebra on $A$ is clearly
determined by the collection $m_i[y]^{\otimes i}, i=1,2,\ldots$.
Notice that the map $\partial$ is necessarily given by a
multiplication by some $x\in R$ so that $\partial(y)=x\cdot 1$. If
$d$ is odd, then $\partial=0$. If $d$ is even and $x$ is not a
zero divisor in $R$, then the (internal) homology of $A$ is simply
$R/x$. For an $E_\infty$ ring spectrum $R$ the structure of the
(homotopy) associative algebra on $R/x$ was investigated in
\cite{EKMM} and \cite{Str}. This parallel topological theory was
our original motivation for introducing the notion of a Moore
algebra.

Let $R^\prime$ be another evenly graded commutative ring and
$f:R\rightarrow R^\prime$ be a ring map. Consider a Moore algebra
$A$ over $R$ specified by the collection $\{m_i[y]^{\otimes i}\in
R\}$. Then the collection $\{f(m_i[y]^{\otimes i})\in R^\prime\}$
will determine a Moore algebra $f_*A$ over $R^\prime$. In other
words the set $\mathcal{S}(d)$ which associates to any evenly
graded commutative ring $R$ the set of Moore algebras over $R$ of
degree $d$ is a functor of $R$.
\begin{theorem}
$(i)$. Let $d$ be even. Then the functor $\mathcal{S}(d)$ is
representable by the polynomial algebra
$R_e=\mathbb{Z}[u_1,u_2,\ldots]$ where $|u_i|=i(d+2)-2$.  More
precisely there exists a Moore algebra $A_e$ over $R_e$ of degree
$d$ such that for any $R$ and any Moore algebra $A$ over $R$ of
degree $d$ there exists a unique ring map $R_o\rightarrow R$ for
which $f_*A_e=A$. The universal Moore algebra $A_e$ is specified
by the formulae  $m_i[y]^{\otimes i}=u_i[1],i=1,2,\ldots$.

$(ii)$. Let $d$ be odd. Then $\mathcal{S}(d)$ is represented by
the polynomial algebra
$R_o=\mathbb{Z}[v_1,v_2,\ldots]\otimes\mathbb{Z}[w_1,w_2,\ldots]$
where $|v_i|=2i(d+2)-d-3$ and $|w_i|=2i(d+2)-2$. More precisely
there exists a Moore algebra $A_o$ over $R_o$ of degree $d$ such
that for any $R$ and any Moore algebra $A$ over $R$ of degree $d$
there exists a unique ring map $R_o\rightarrow R$ for which
$f_*A_o=A$. The universal Moore algebra $A_o$ is specified by the
formulae $m_{2i-1}[y]^{\otimes 2i-1}=0$ and $m_{2i}[y]^{\otimes
2i}= v_i[y]+w_i[1], i=1,2,\ldots$.
\end{theorem}
\begin{proof} In both cases $(i)$ and $(ii)$ the universality is
obvious and we only need to prove that $m\circ m=0$. Note that
apriori the latter equation could impose nontrivial relations on
the generators $u_i, v_i$ and $w_i$; the theorem effectively
states that no such relations except commutativity are in fact
present.

The equality  $m\circ m=0$ could be checked directly using the
composition formula (\ref{big}). However this path is rather
long-winded and unenlightening and we will choose the approach via
the cobar-construction. So consider the algebra $(T\Sigma
A)^\ast=R\langle\langle\tau,t\rangle\rangle$ where $\tau$ and $t$
are dual to $[1]$ and $[y]$ respectively so $|\tau|=-1,|t|=-d-2$.
Then the coderivation $m$ of $T\Sigma A$ determines the continuous
derivation $m^\ast$ of $(T\Sigma A)^\ast$. Routine inspection
shows that in the case $(i)$
$$m^\ast=\sum_{i=1}^{\infty}u_it^i\partial_\tau+ad\tau-\tau^2\partial_{\tau}$$
whereas in the case $(ii)$ we have
$$m^\ast=\sum_{i=1}^{\infty}v_it^{2i}\partial_t+\sum_{i=1}^{\infty}
w_it^{2i}\partial_\tau-\tau^2\partial_\tau+ad\tau.$$ For $(i)$ we
compute $$(m^\ast\circ m^\ast)(t)=m^\ast([\tau,t])=
[\sum_{i=1}^{\infty}u_it^i+\tau^2,t]-[\tau,[\tau,t]].$$ The
elements $t$ and $\sum_{i=1}^{\infty}u_it^i$ are both even and
therefore $[\sum_{i=1}^{\infty}u_it^i,t]=0$. Clearly
$[\tau^2,t]-[\tau,[\tau,t]]=0$ This implies that  $(m^\ast\circ
m^\ast)(t)=0$. Next, \begin{eqnarray}\nonumber(m^\ast\circ
m^\ast)(\tau)=m^\ast(\sum_{i=1}^{\infty}u_it^i+\tau^2)\\
\nonumber=
[\tau,\sum_{i=1}^{\infty}u_it^i]+m^\ast(\tau^2)=[\tau,\sum_{i=1}^{\infty}u_it^i]+
\sum_{i=1}^{\infty}u_it^i\partial_{\tau}(\tau^2)=0\end{eqnarray}
It proves that $m^\ast\circ m^\ast=0$. Similarly for $(ii)$ we
have
\begin{eqnarray}\nonumber(m^\ast\circ
m^\ast)(t)=m^\ast(\sum_{i=1}^{\infty}v_it^{2i}+[\tau,t])\\
\nonumber
=(\sum_{i=1}^{\infty}v_it^{2i}\partial_{t})(\sum_{i=1}^{\infty}v_it^{2i})-
[\tau,\sum_{i=1}^{\infty}v_it^{2i}]+[\tau,\sum_{i=1}^{\infty}v_it^{2i}]+
[\tau^2,t]-[\tau,[\tau,t]].\end{eqnarray}
 Since now the element
$t$ is odd the derivation~
$\sum_{i=1}^{\infty}v_it^{2i}\partial_{t}$~ is also odd while
$\sum_{i=1}^{\infty}v_it^{2i}$~ is even  and it follows that
$(\sum_{i=1}^{\infty}v_it^{2i}\partial_{t})(\sum_{i=1}^{\infty}v_it^{2i})=0$.
Just as before we have $[\tau^2,t]-[\tau,[\tau,t]]=0$. Therefore
$(m^\ast\circ m^\ast)(t)=0$. Further
\begin{eqnarray}\nonumber(m^\ast\circ
m^\ast)(\tau)=m^\ast(\sum_{i=1}^{\infty}w_it^{2i}+\tau^2)\\
\nonumber=
(\sum_{i=1}^{\infty}v_it^{2i}\partial_{t})(\sum_{i=1}^{\infty}w_it^{2i})+
(\sum_{i=1}^{\infty}w_it^{2i}\partial_{\tau})(\tau^2)+
[\tau,\sum_{i=1}^{\infty}w_it^{2i})].\end{eqnarray} Arguing as
before we see that the first term in the last expression is zero
whereas the second and third cancel each other out. Therefore
$(m^\ast\circ m^\ast)(\tau)=0$ and we are done.
\end{proof}
\begin{rem}For an odd $d$ the differential on the underlying
complex $A=\{\Sigma^dR\rightarrow R\}$ is zero and therefore its
homology is fixed. For $d$ even the differential is given by
multiplication with $u_1=x\in R_e$. The element $u_1$ plays a
special role among $u_i$'s fixing the homology of the Moore
algebra. We will be interested mostly in the case when $x$ is a
nonzero divisor in $R$ in which case $H_\ast(A)=R/x$.\end{rem}
\begin{rem} The universal odd Moore algebra has an ideal generated
by the element $y$. This is a nonunital $A_\infty$-algebra over
$R_o$ such that $m_{2i}[y]^{2i}=v_i[y]$ and $m_{2i-1}=0$. This
$A_\infty$-algebra was introduced in the early nineties by
M.Kontsevich.  It turns out to be related to Morita-Miller-Mumford
classes in the cohomology of moduli spaces of algebraic curves,
cf. \cite{Kon}. It would be interesting to understand whether our
more general constructions can yield new information about
cohomologies of these moduli spaces.\end{rem} We see, therefore,
that an arbitrary even Moore algebra $A$ over a ring $R$ is
specified by the collection $\{u^A_i\}\in R$ where $u^A_i$ is the
image of $u_i$ under the classifying map $R_e\rightarrow R$. In
that case the $A_\infty$ structure on $A$ is the following
derivation $m_A^\ast$ of the algebra $(T\Sigma
A)^\ast=R\langle\langle\tau,t\rangle\rangle$:
$$m^\ast_A=\sum_{i=1}^{\infty}u_i^At^i\partial_\tau+ad\tau-\tau^2\partial_{\tau}.$$
Similarly an odd Moore algebra over $R$ is determined by the
collection $\{v_i^A,w_i^A\in R\}$, the images of $v_i$ and $w_i$
under the classifying map $R_o\rightarrow R$. The $A_\infty$
structure on $A$ is the following derivation $m_A^\ast$ of the
algebra $(T\Sigma A)^\ast$:
$$m^\ast_A=\sum_{i=1}^{\infty}v^A_it^{2i}\partial_t+\sum_{i=1}^{\infty}
w^A_it^{2i}\partial_\tau-\tau^2\partial_\tau+ad\tau.$$ We see that
an even (odd) Moore algebra is completely characterized by a power
series $u^A(t):=\sum_{i=1}^{\infty}u_i^At^i$ (by a pair of power
series
$(v^A(t),w^A(t)):=(\sum_{i=1}^{\infty}v^A_it^{2i},\sum_{i=1}^{\infty}w^A_it^{2i})$
respectively). We will call these power series \emph{
characteristic power series} for corresponding Moore algebras.
\begin{rem} The notion of a characteristic power series is similar to
that of a formal group law. Further the universal (even or odd)
Moore algebra is analogous to the universal formal group law over
the Lazard ring (which is also a polynomial ring in infinitely
many variables). The Moore algebras corresponding to different
points of $R_e$ or $R_o$ could still be weakly equivalent (note
that this is exactly what happens also for formal groups).
Moreover we have certain infinite-dimensional Lie groups acting on
$R_o$ and $R_e$ whose orbits correspond to weakly equivalent Moore
algebras. These actions are far from being free which means that
there are moduli stacks rather than moduli spaces of Moore
algebras.  We will see in the next section that in the even case
the corresponding group is just the group of formal power series
in one variable with vanishing constant term. This again forces
one to think of the analogy with formal groups.\end{rem}
 \section{Classification problem}
It is an interesting and nontrivial problem to classify Moore
algebras over a given ring up to a (unital) weak equivalence. In
this paper we will consider only the even case. An even Moore
algebra $A$ of degree $d$ has the characteristic series
\begin{equation}\label{car}
u^A(t)=u(t)=\sum_{i=1}^{\infty}u_it^i\end{equation} Here
$|t|=-(d+2)$ and $|u_i|=i(d+2)-2$ from which we conclude that
$|u(t)|=-2$. Conversely any such power series determines an even
Moore algebra. It is easy to see that if
$f=(f_1,f_2,\ldots):BA\rightarrow BC$ is a weak equivalence
between two even Moore algebras then $f_1$ is an isomorphism so
$f$ is invertible (even though $A$ and $C$ need not be minimal).
It follows that the set of weak equivalence classes of even Moore
algebras coincides with the set of orbits of the group
$\overline{Aut}(R\langle\langle\tau,t\rangle\rangle)$ on the set
of (unital) $A_\infty$ structures on $A$ which could be identified
with the set of power series (\ref{car}). According to Remark
\ref{act1} the group
$\overline{Aut}(R\langle\langle\tau,t\rangle\rangle)$ is the group
of formal power series $f(t)=f_1t+f_2t^2+\ldots$, where $f_1$ is
invertible and the group operation is composition. (Notice that
the condition that $\overline{Aut}(R\langle\langle\tau,{\bf
t}\rangle\rangle)$ consists of morphisms of zero degree imposes
some restrictions on $f_i$, namely $|f_i|=(i-1)(d+2)$.) The action
is given by substitution of power series. To summarize we have the
following
\begin{theorem} The set of equivalence classes of even Moore
algebras over $R$ is in $1-1$ correspondence with the set of
orbits of the group $\overline{Aut}(R\langle\langle\tau,{\bf
t}\rangle\rangle)$ acting on the set of formal power series with
coefficients in $R$ of degree $-2$ with vanishing constant term.
The action of the group element $f(t)$ on the power series $u(t)$
is given by the formula $u(t)\rightarrow u(f(t))$.
\end{theorem}
\begin{rem} Suppose that the ring
$R$ is $2$-periodic, i.e. it possesses an invertible element $v$
of degree $2$. In that case the group
$\overline{Aut}(R\langle\langle\tau,{\bf t}\rangle\rangle)$ is
isomorphic to the group of formal power series with coefficients
in $R_0$, the zeroth component of $R$ and having vanishing
constant term. The set of characteristic series becomes the set of
all power series with coefficients in $R_0$ without constant term
and the action is given by substitution as above. In other words
there are no degree restrictions on the coefficients of power
series.\end{rem}
 Just as for formal
groups it seems hopeless to try to make the classification  over
an arbitrary ring. The restriction that we place on $R$ is that we
assume that $R$ is either a (graded) field or a (graded) complete
discrete valuation ring. We refer the reader to the book
\cite{Ser} by Serre for an account on discrete valuation rings. In
this book the ungraded rings are treated but passage to the graded
case is automatic.

Next we introduce the notion of the height of a formal power
series which will be one of the invariants of the associated Moore
algebra.
\begin{defi}Let $u(t)=\sum_{i=1}^{\infty}u_it^i$ be a formal power
series without a constant term. Then we say that $u(t)$ has height
$n$ if $u_n$ is the first nonzero coefficient of $u(t)$.  The
height of the characteristic series of an even Moore algebra $A$
is called the height of $A$.
\end{defi}
\begin{prop}\label{field}Let $R$ be a graded field of characteristic zero, $
A$ and $C$ be two even Moore algebras over $R$  with
characteristic series $u^A(t)=\sum_{i=1}^{\infty}u^A_it^i$ and
$u^C(t)=\sum_{i=1}^{\infty}u^C_it^i$. Then $A$ is weakly
equivalent to $C$ if and only if  $n=$height$(A)=$height$(C)$ and
$r^nu^A_n=u_n^C$ for some $r\in R_0$. Thus the equivalence class
of an even Moore algebra of degree $d$ is determined by a pair
$(n,r)$ where $n$ is the height and $r\in
R_0^\times/{R_0^\times}^n$ is an element in $R_0^\times$ modulo
the subgroup of $n$th powers.\end{prop}
\begin{proof}Let $A$ have height $n$.
Then $u^A(t)=\sum_{i=n}^{\infty}u^A_it^i$. First we prove that
there exists a power series $h(t)$ such that $u^A(h(t))=u^A_nt^n$.
Let $k_1\in \mathbb{Z}$ the smallest integer for which
$u^A_{k_1}\neq 0$ and  $k_1>n$. If no such integer exists then
$u^A(t)$ is already in the desired form $u^A(t)=u^A_nt^n$.
Otherwise consider the polynomial
$$h_1(t)=t-\frac{t^{k_1-(n-1)}u^A_{k_1}}{nu^A_n}.$$ Then by Taylor's
formula $$u^A\circ
h_1(t)=u^A(t)-[\frac{d}{dt}u^A(t)]\frac{t^{k_1-(n-1)}u^A_{k_1}}{nu^A_n}\mod(t^{k_1+1})
=u_n^At^n \mod(t^{k_1+1}).$$ Now let $k_2>n$ be the smallest
integer for which the coefficient at $t^{k_2}$ in $u^A\circ
h_1(t)$ is nonzero. Clearly $k_2>k_1$. Then just as before we
could find $h_2(t)$ for which $u^A\circ h_1\circ
h_2(t)=u^A_nt^n\mod(t^{k_2})$. Continuing this process we
construct the sequence $h_1,h_2,\ldots$ of polynomials of the form
$h_i=t+a_it^{k_i}$ for some $a_i\in R$. The sequence $\{h_1\circ
h_2\circ\ldots\circ h_i\}$ clearly converges in the $t$-adic
topology and denoting its limit by $h(t)$ we obtain
$u^A(h(t))=u^A_nt^n$.

In other words we proved that $A$ is weakly equivalent to the
Moore algebra $A^\prime$ having the characteristic series
$u^{A^\prime}(t)=u^A_nt^n$. Similarly $C$ is equivalent to a Moore
algebra $C^\prime$ with characteristic series
$u^{C^\prime}(t)=u^C_mt^m$. The rest is clear: $A^\prime $ and
$C^\prime$ are equivalent if and only if \begin{enumerate} \item
$n=m$ and \item
$u^{A^\prime}(rt)=r^nu^A_nt^n=u^{C^\prime}(t)=u_n^Ct^n$ for some
$r\in R_0^\times$\end{enumerate} which means that
$r^nu^A_n=u_n^C$.
\end{proof}
\begin{rem} The assumption that $R$ has characteristic $0$ could be replaced with
the assumption that height($A$)=height($C$) does not divide
$char(R)$. The proof is the same verbatim.\end{rem}
\begin{rem} Notice that the statement of Theorem \ref{field} is vacuous in the
case of height $1$. Indeed, an even Moore algebra having the
characteristic series $\sum_{i=1}^{\infty}u_it^i$ with $u_1$
invertible is trivial since the underlying complex
$\Sigma^dR\rightarrow R$ is contractible. To get a nontrivial even
Moore algebra of height $1$ we have to have nonzero noninvertible
(homogeneous) elements in the ground ring $R$. This is the
situation that arises in the study of $MU$-modules and
$MU$-algebras in stable homotopy theory, cf. \cite{Str} and
\cite{Laz}. In order to obtain reasonable classification results
we need to impose certain conditions on $R$.\end{rem}
\begin{defi}Let $R$ be a (graded) complete discrete valuation ring with
uniformizer $\pi$ and $u(t)=\sum_{i=1}^{\infty}u_it^i$ is a power
series with coefficients in $R$. We will call $u(t)$ {\it trivial}
if $u(t)=\pi t$ and {\it canonical} if there exists an $n$ for
which \begin{enumerate}\item $u_1=\pi$; \item $\pi$ divides
$u_2,u_3,\ldots,u_{n-1}$; \item $u_n$ is invertible; \item
$u_{n+1}=u_{n+2}=\ldots=0$.\end{enumerate}\end{defi}
\begin{rem}A canonical power series $u(t)$ could be defined
equivalently as $tP(t)$ where $P(t)$ is an Eisenstein
polynomial.\end{rem}
\begin{prop}\label{val}Let $R$ be a graded complete discrete valuation ring with
residue field of characteristic $0$ and uniformizer $\pi$. Let $A$
be an even Moore algebra having characteristic series
$u^A(t)=\sum_{i=1}^{\infty}u_i^At^i$ where $u_1^A=r\pi$, $r\in
R^\times$. Then $A$ is weakly equivalent to the algebra having
either trivial or canonical characteristic series. Moreover two
Moore algebras having canonical characteristic series are weakly
equivalent if and only if these series coincide.
\end{prop}
\begin{rem}Recall, that a complete discrete valuation ring $R$ with
residue field $R/\pi$ of characteristic zero is isomorphic to the
formal power series ring $R/\pi[[T]]$, in particular $R$ is a
vector space over $R/\pi$.\end{rem}
\begin{proof}In the interests of readability we suppress the
superscript $A$ and will write $u(t)=\sum_{i=1}^{\infty}u_it^i$
for $u^A(t)=\sum_{i=1}^{\infty}u_i^At^i$. First we could assume
that $u_1=\pi$ or else use the substitution $u(t)\rightarrow
u(r^{-1}t)$ to reduce $u(t)$ to the desired form. Next suppose
that all coefficients of $u(t)$ are divisible by $\pi$. Then
clearly using substitutions $u(t)\rightarrow u(t+a_nt^n)$ for
suitable $n$ and $a_n$ we could eliminate these coefficients one
by one and reduce $u(t)$ to the trivial form. Now assume that not
all $u_i$ are divisible by $\pi$ and denote by $u_k$ the first
such. In other words $u(t)=\sum_{i=1}^{k}u_it^i\mod(t^{k+1})$
where $u_i=0~\mod(\pi)$ for $i=1,2,\ldots,k-1$ and $u_k$ is
invertible in $R$.

If for $i>k$ the coefficients $u_i$ are zero we are done since
$u(t)$ is already in the canonical form. If not let $l(u)$ be the
maximal integer $l$ for which $u_i=0~\mod(\pi^l)$,
$i=k+1,k+2,\ldots$. Notice that $l(u)$ could be zero. Our first
step would be to find an appropriate substitution $u(t)\rightarrow
u(h^1(t))$ such that $l(u(h))>l(u)$. Set $$s_1:=min\{i:i>k,
\pi^{l(u)} \mbox{~divides~} u_i \mbox{~but~}
\pi^{l(u)+1}\mbox{~does not divide~}u_i\}.$$ Then we have
$$u(t)=\sum_{i=1}^{k}u_it^i+u_{s_1}t^{s_1}\mod(t^k\pi^{1(u)+1}+t^{s_1+1}\pi^{l(u)}).$$
Let $h_1(t):=t-\frac{u_{s_1}}{ku_k}t^{s_1-(k-1)}$ (recall that
$u_k$ is invertible). Then Taylor's formula implies that $$u\circ
h _1(t)=\sum_{i=1}^{k}v_it^i
\mod(t^k\pi^{1(u)+1}+t^{s_1+1}\pi^{l(u)}) $$ where
$v(t):=\sum_{i=1}^{k}v_it^i$ is a canonical polynomial. Notice
that $l(v)\geq l(u)$. If $l(v)>l(u)$ then our first step is
completed. Assuming that $l(v)=l(u)$ set $$s_2:=min\{i:i>k,
\pi^{l(v)} \mbox{~divides~} v_i \mbox{~but~}
\pi^{l(v)+1}\mbox{~does not divide~}v_i\}.$$ Observe that
$l(u)=l(v)$ implies $s_2>s_1$. It follows that
$$v(t)=\sum_{i=1}^{k}v_it^i+u_{s_2}t^{s_2}
\mod(t^k\pi^{1(u)+1}+t^{s_2+1}\pi^{l(u)}).$$ Then just as before
set $h_2(t):=t-\frac{v_{s_2}}{kv_k}t^{s_2-(k-1)}$ and consider the
series $w(t):=v(h_2(t))=u\circ h_1\circ h_2$. Continuing in this
way we construct a sequence of power series $\{u\circ h_1\circ
h_2\circ\ldots\circ h_n\}$. This is clearly a Cauchy sequence in
the $t$-adic topology and converges (or perhaps stops at a finite
stage) to a power series $u^1(t)$ having the property that
$l(u^1)>l(u)$. Notice that $u^1(t)=u(h^1(t))$ where $h(t)=h_1\circ
h_2\circ\ldots$. Moreover we have $h^1(t)=t~\mod(\pi^{l(u)})$.
This completes the first step. It is clear how to proceed further.
Repeating the above procedure we find  power series $h^2(t)$ and
$u^2(t):=u^1(h^2(t))$ such that $h^1(t)=t~\mod(\pi^{l(u^1)})$ and
$l(u^2)>l(u^1)$. The sequence $\{h^1\circ h^2\circ\ldots\circ
h^n\}$ is a Cauchy sequence in the $\pi$-adic topology and
converges to $h(t)$. Then $u(h(t))$ is a canonical polynomial.

We still need to prove that two even Moore algebras having
different canonical polynomials cannot be equivalent. In other
words we have to show that if $u(t)=\sum_{i=1}^{n}u_it^i$ and
$v(t)=\sum_{i=1}^{m}v_it^i$ are two canonical polynomials and
$h(t)=\sum_{i=1}^{\infty}h_it^i$ is such that
\begin{equation}\label{sub}u(h(t))=v(t)\end{equation} then $u(t)=v(t)$.
Indeed since $u_1=v_1=\pi$ the equality (\ref{sub}) implies that
$h_1=1$ so we have $h(t)=t+\tilde{h}(t)$ where
$\tilde{h}(t)=0\mod(t^2)$. Further from (\ref{sub}) we obtain
$u_n(h(t))^n=v_mt^m\mod(\pi)$. It follows that $m=n$, $u_n=v_m\mod
(\pi)$ and $(h(t))^n=t^n\mod(\pi)$. Since the residue field
$R/\pi$ has characteristic $0$ the last equality implies
$h(t)=t~\mod(\pi)$ or equivalently $\tilde{h}(t)=0~\mod(\pi)$. Now
suppose that $\tilde{h}\neq 0$ and  let $k$ be the unique integer
for which $\pi^k$ divides $\tilde{h}(t)$ and $\pi^{k+1}$ does not
divide $\tilde{h}(t)$. By Taylor's formula
$$u(h(t))=u(t+\tilde{h}(t))=
u(t)+u^\prime(t)\tilde{h}(t)~\mod(\pi^{k+1}).$$ Since $u_k\neq
0~\mod(\pi)$ the series $u^\prime(t)\tilde{h}(t)~\mod(\pi^{k+1})$
will necessarily have nonzero terms of order $>k$ in $t$. This is
a contradiction with our assumption that $u(h(t))=v(t)$ is a
canonical polynomial of degree $k$. With this our proposition is
proved.
\end{proof}
\begin{rem}One  would naturally like to know to whether
Proposition \ref{val} remains true if the residue field $R/\pi$
has characteristic $p$. Suppose that is the case and let $A$ be an
even Moore algebra having characteristic series
$u^A(t)=\sum_{i=1}^{\infty}u_i^At^i$ with $u_1=r\pi, r\in
R^\times$. If all coefficients $u_i, i=2,3,\ldots$ of $u^A(t)$ are
divisible by $\pi$ then exactly as in the proof of Proposition
\ref{val} one shows that $A$ is equivalent to a Moore algebra
having characteristic series $u(t)=t$. If not, let $u^A_k$ be the
first invertible coefficient in $u(t)$. If $p$ does not divide $k$
then the proof of Proposition \ref{val} carries over verbatim to
show that $u^A(t)$ can be reduced to a canonical form and this
canonical form is unique. Therefore in this case we have an exact
analogue of Proposition \ref{val}.  If $p$ does divide $k$ the
classification seems to be much more subtle. \end{rem}

\section{Cohomology of Moore algebras}
In this section we compute  Hochschild cohomology of Moore
algebras of even degree subject to the condition that the first
coefficient of its characteristic series is a nonzero divisor  and
discuss their connection with totally ramified extensions of local
fields. It would be interesting to calculate Hochschild cohomology
for {\it odd} Moore algebras. In principle the same method should
apply, however in order to get a sensible answer one has to place
some restrictions on characteristic series and it is not
immediately clear what these restrictions should be.

\begin{prop}\label{der} Let $A$ be the Moore algebra of even degree  over $R$ with
characteristic series $u(t)=u^A(t)$. Let us assume that the
coefficient $u_1=u_1^A$ of $u(t)$ is not a zero divisor in $R$.
Then there is an isomorphism of $R$-modules $HH^\ast(A,A)\cong
R[[t]]/(u^\prime(t))$ where $u^\prime(t)$ denotes the derivative
of the power series $u(t)$.\end{prop}
\begin{proof} We will compute $HH^\ast(A,A)$ as the homology of the operator
$ [?,m^\ast]$ on the space of normalized (continuous) derivations
of $T\Sigma A^\ast=R\langle\langle\tau,t\rangle\rangle$. Recall
that $m^\ast=u(t)\partial_\tau+ad\tau-\tau^2\partial_{\tau}.$  Let
$\xi$ be a normalized derivation so
$\xi=A(t)\partial_{\tau}+B(t)\partial_t$. Then completely
automatic calculations show that
$[\xi,m^\ast]=u^\prime(t)B(t)\partial_\tau$. The condition that
$u_1$ is not a zero divisor in $R$ implies that $u^\prime(t)$ is
not a zero divisor in $R[[t]]$. Therefore the kernel of the
operator $[?,m^\ast]$ consists of derivations of the form
$A(t)\partial_\tau$ whereas its image is precisely the derivations
of the form $u^\prime(t)B(t)\partial_\tau$ and we are done.
\end{proof}

Let us now look at the Hochschild cohomology of an even Moore
algebra $A$ from the point of view of the associated dga
$\tilde{A}$. Without loss of generality we may suppose that
$\tilde{A}$ is a cell $R$-module.  Again, our standing assumption
is that $u^A_1$ is not a zero divisor in $R$ so the internal
homology of $A$ is $R/u_1^A$. We have the classical Hochschild
(bi)complex ${C}^\ast(\tilde{A},\tilde{A})$:
\begin{equation}\label{Ho}{\tilde{A}}\rightarrow Hom({\tilde{A}}
,\tilde{A})\rightarrow\ldots\rightarrow Hom({\tilde{A}}^{\otimes n
},\tilde{A})\rightarrow \ldots\end{equation} Associated to this
bicomplex is the spectral sequence with the $E^1$-term $E^1_{\ast
n }=H_\ast(Hom(\tilde{A}^{\otimes n},\tilde{A}))$. Since
$\tilde{A}$ is weakly equivalent to a finite cell $R$-module the
natural map $Hom(\tilde{A},R)\otimes \tilde{A}\rightarrow
Hom(\tilde{A},\tilde{A})$ is a homology isomorphism. We have the
following sequence of homology isomorphisms:
\begin{eqnarray}\nonumber Hom(\tilde{A}^{\otimes 2},\tilde{A})\simeq
Hom(\tilde{A},R)\otimes Hom(\tilde{A},\tilde{A})\\ \nonumber\simeq
Hom(\tilde{A},R) \otimes\tilde{A}\otimes_{\tilde{A}}
Hom(\tilde{A},\tilde{A})\simeq
Hom(\tilde{A},\tilde{A})\otimes_{\tilde{A}}Hom(\tilde{A},\tilde{A}).\end{eqnarray}
More generally we have the following homology isomorphism
$$Hom(\tilde{A}^{\otimes n},\tilde{A})\simeq
Hom(\tilde{A},\tilde{A})\otimes_{\tilde{A}}Hom(\tilde{A},\tilde{A})
\otimes_{\tilde{A}}\ldots
\otimes_{\tilde{A}}Hom(\tilde{A},\tilde{A})\mbox{ ($n$ times}).$$
Further a straightforward computation shows
$$H_\ast
Hom(\tilde{A},\tilde{A})=Ext^\ast_R(R/u_1^A,R/u_1^A)=\Lambda_{R/u_1^A}(z)$$
where $\Lambda_{R/u_1^A}(z)$ denotes the exterior algebra over
$R/u_1^A$ on one generator $z$ of degree $-|u_1^A|-1=-d-1$. So the
$E_1$-term of our spectral sequence has the form
$$H_\ast(A)=R/u_1^A\rightarrow
\Lambda_{R/u_1^A}(z)\rightarrow\ldots\rightarrow
(\Lambda_{R/u_1^A}(z))^{\otimes n}\rightarrow \ldots$$ This is the
usual cobar complex for the Hopf algebra $H_\ast
Hom(\tilde{A},\tilde{A})=\Lambda_{R/u_1^A}(z)$ and its homology is
\begin{equation}\label{rf}E_2^{\ast\ast}=
Ext^{\ast\ast}_{\Lambda_{R/u_1^A}(z)}(R/u_1^A,R/u_1^A)=R/u_1^A[[t]]
\end{equation} where $t$ has degree $-d-2$. For dimensional reasons
$E^2=E^3=\ldots=E^\infty$.
 Next
notice that the the spectral sequence $E_1^{\ast\ast}$ is
multiplicative via the pairing
$$Hom({\tilde{A}}^i,\tilde{A})\otimes
Hom({\tilde{A}}^j,\tilde{A})\rightarrow
Hom({\tilde{A}}^{i+j},\tilde{A}\otimes\tilde{A})\rightarrow
Hom({\tilde{A}}^{i+j},\tilde{A})$$ where the second map is induced
by the multiplication $\tilde{A}\otimes \tilde{A}\rightarrow
\tilde{A}$. This pairing turns $E_1^{\ast\ast}$ into a graded ring
and it follows that (\ref{rf}) is in fact an isomorphism of rings.

So we proved the following
\begin{prop}\label{hoch}The Hochschild cohomology ring of an even Moore
algebra of degree $d$ over $R$ with $u_1^A\in R$ a nonzero divisor
is a complete filtered ring whose associated graded ring is a
formal power series algebra over $R/u_1^A$ on one generator in
degree $-d-2$.
\end{prop}
\begin{rem}The arguments above would be considerably simpler if we
knew that the spectral sequence of Proposition \ref{ss} were
multiplicative. Unfortunately this is not known yet. \end{rem}

Now let $R$ be a (graded) discrete valuation ring with the
uniformizer $\pi=u_1$ and consider the unit map $f:R\rightarrow
H^\ast(A,A)$ where $A$ is as in Proposition \ref{hoch}. Since the
filtered ring $H^\ast(A,A)$ has a formal power series ring over
the field $R/(\pi)$ for its associated graded ring we conclude
that $H^\ast(A,A)$ has no zero divisors. Therefore the map $f$ is
either an injection or its kernel is the maximal ideal
$(\pi)\subset R$. Furthermore the ring $H^\ast(A,A)$ is itself a
graded discrete valuation ring. Using Proposition \ref{der} we see
that the kernel of $f$ is $(\pi)$ if and only if
$u^\prime(t)=0~\mod(\pi)$. The last equality is equivalent to
$u(t)=0~\mod(\pi)$ if $char(R/\pi)=0$ or to $u(t)=v(t^p)\mod
(\pi)$ for some power series $v(t)$ if $char(R/\pi)=p$.

Now suppose that $f$ is injective. Then Proposition  \ref{hoch}
implies that $u^\prime(t)$ is not divisible by $\pi$ which means
that there exists $n\in\mathbb{Z}$ for which $nu_n$ is invertible
in $R$. Consider the smallest such $n$; it obviously equals the
height of the series $u(t)$ reduced $\mod(\pi)$. By the
Weierstrass Preparation Theorem the ring $H^\ast(A,A)$ is free of
rank $n$ over $R$. Therefore we obtain the following
\begin{cor}Let $R$ be a graded discrete valuation ring
with uniformizer $\pi$ and residue field $R/\pi$ of characteristic
$p$. Let $A$ be an even Moore algebra over $R$ with characteristic
series $u(t)=\sum_{i=1}^{\infty}u_i^At^i$ where $u_1^A=\pi$. Then

 $(i)$ The ring $H^\ast(A,A)$ is either
an $R/\pi$-algebra or a totally ramified extension of $R$ whose
ramification index equals the height of the series
$u(t)~\mod(\pi)$.

 $(ii)$The ring $H^\ast(A,A)$ is an
$R/\pi$-algebra if and only if $u(t)=v(t^p)\mod(\pi)$ for some
polynomial $v(t)$.\end{cor}
 \begin{rem}Varying $u(t)$ we could get ramified extensions of arbitrary
index that is coprime to $p$. In particular if $charR/\pi\neq 2$
the inclusion $f:R\hookrightarrow H^\ast(A,A)$ could be an
isomorphism. $A_\infty$-algebras having the property that $f$ is
an isomorphism are analogous to central separable algebras which
were studied extensively in ring theory, cf. \cite{AG} and we hope
to return to them in the future.\end{rem}

We conclude this section with a few simple examples illustrating
our results. Let $R=\hat{\mathbb{Z}}_p[v,v^{-1}]$ where
$\mathbb{Z}_p$ is the ring of $p$-adic integers for $p\neq 2$ and
$v$ is a formal Laurent variable of degree $2$. Consider two
$A_\infty$ structures $m^1$ and ${m}^\infty$ on the complex
$A=\{R\stackrel{p}{\rightarrow} R\}$. These will in fact be
differential graded algebra structures, i.e.
$m^1_i={m}^\infty_i=0$ for $i>2$. Namely, set $m_2^\infty[y|y]=0$
and ${m}_2^1[y|y]=v[1]$. In other words  $(A,m^\infty)$ is just
the exterior algebra on $y$ in degree $1$ with differential $dy=p$
while $(A,{m}^1)$ is the dga generated by $y$ with the same
differential $dy=p$ but with the relation $y^2=v$. Then Hochschild
cohomology of $(A,m^\infty)$ is just the algebra
$R/p[[t]]=\mathbb{F}_p[v,v^{-1}][[t]]$ while Hochschild cohomology
of $(A,{m}^1)$ is the ring $R$ itself (the ramification index
equals $1$ in this case). The notations $m^1$ and $m^\infty$
suggest that there are also $m^n$'s for finite $n$. These indeed
exist and could be obtained by setting $m^n_i[y]^{\otimes i}=0$
for $i\neq n$ and $m^n_n[y]^{\otimes n}=v^n[1]$. The Hochschild
cohomology of $(A,{m}^n)$ realizes a totally ramified extension of
the $p$-adic integers of index $n$ which is not divisible by $p$.
\\
\\{\bf Concluding remarks.} It should be noted that our present
approach to the moduli problem is rather ad hoc and it would be
valuable to consider it from the more general point of view. Here
we mention the (still unpublished) preprint of M. Schlessinger and
J. Stasheff \cite{SS} where this program is carried out for
rational homotopy types. These authors effectively study the {\it
commutative} $A_\infty$ structures on a complex with fixed
homology ring $H$ over the field of rationals. They replace $H$
with its multiplicative resolution $\Lambda Z$, and consider the
graded Lie algebra $Der(\Lambda Z)$  of derivations of $H$. Then
it turns out that the moduli space under consideration is
represented by the standard construction $A(Der(\Lambda Z))$ which
computes homology of the Lie algebra $Der(\Lambda Z)$. This simple
and elegant approach is very appealing and we feel that it is
possible to extend it in the context of $A_\infty$-algebras. The
role of $Der(\Lambda Z)$ should be played by the Hochschild
complex $C^\ast(H,H)$.

It is now clear that the set of homotopy types of dga's with a
fixed homology algebra is only $\pi_0$ of the `true' moduli space.
The other invariants are picked up by monoids of homotopy
self-equivalences corresponding to different path components of
the moduli space. This point of view is developed in \cite{BDG}.
However in this context the problem of computing $\pi_0$ differs
sharply from that of computing the higher homotopy groups. Indeed,
in \cite{Laz3} we showed that, essentially, higher homotopy groups
of mapping spaces could be reduced to (a version of) Hochschild
cohomology. In fact in the cited reference the result is obtained
for $S$-algebras but the arguments are still valid for dga's.

 Also relevant to this problem is the recent
paper by V. Hinich \cite{Hin} where homotopy invariant deformation
theory was constructed in the abstract setting of an algebra over
an operad. However the approach in the cited reference is
restricted by working in characteristic $0$ and considering
connected algebras only.

Another problem is to extend our results to the category of
$R$-algebras in the sense of \cite{EKMM}. In the simplest case,
which is already highly nontrivial, one is asked to classify the
structures of $KU$-algebra structures on $KU/p$. Here $KU$ is the
spectrum of topological $K$-theory which is known to be a
commutative $S$-algebra. The related problem is to compute
$THH(KU/p,KU/p)$, the topological Hochschild cohomology of $KU/p$.
We saw that in the algebraic case we obtain tamely ramified
extensions of the $p$-adics. Perhaps in the topological case one
could get wildly ramified extensions?

\end{document}